\documentclass[notitlepage,10pt,NFSS]{article}
\usepackage{amsfonts}
\usepackage{amsmath}
\usepackage{amssymb}
\usepackage{theorem}

\theoremstyle{plain}

\newtheorem{Th}{Theorem}[subsection]

\newtheorem{Lem}{Lemma}[subsection]

\newtheorem{Prop}{Proposition}[subsection]

\newtheorem{Cor}{Corollary}[subsection]

{\theorembodyfont{\rmfamily}\newtheorem{Exam}{Example}[subsection]}

{\theorembodyfont{\rmfamily}\newtheorem{Rem}{Remark}[subsection]}

{\theorembodyfont{\rmfamily}\newtheorem{Def}{Definition}[subsection]}
\mathchardef\d="0064
\def\dis{\displaystyle}
\allowdisplaybreaks

\begin{document}

\font\tenrm=cmbx10 at 10pt

\begin{titlepage}
BJGA, Vol. 7, No. 2, 2002, pp. 87-136.

\vskip 0.3cm
\begin{center}
{\LARGE  {\bf{Weierstrass-type Representation of Weakly Regular
Pseudospherical Surfaces in Euclidean Space}}}
\end{center}

\vskip 0.3cm
\begin{center}
{\bf {\Large Magdalena Toda}}
\end{center}

\vskip 2.0cm

\abstract{In [To], the author presented a method of constructing
all weakly regular pseudospherical surfaces corresponding to given
Weierstrass-type data. While the construction itself will appear
later as a separate publication, this report contains a complete
and detailed description of the Weierstrass representation for
weakly regular surfaces with $K=-1$, in terms of moving frames and
loop groups.}

\vskip 2.0cm

{\bf Mathematics Subject Classification}: 53A10, 53C42, 58D10,
49Q05.

\vskip 2.0 cm

{\bf Keywords}: moving frames, pseudospherical surfaces, loop
groups.

\end{titlepage}

\section{Moving Frames of Surfaces in $E^3$}

This is a general introduction to the concept of a moving frame for a surface
in $E^3$,
in the spirit of [Ei] and [Ch, Te].

In the real Euclidean three-space $E^3$ endowed with the inner product
$\langle\cdot,\cdot\rangle$, a {\it frame} is an ordered quadruple
$F=\{x,e_1,e_2,e_3\}$, where $x\in E^3$ and $e_1$, $e_2$, $e_3$ are
orthonormal vectors of positive orientation, i.e., $e_3=e_1\times e_2$. Let
$\cal F$
denote the set of all frames. We will mostly be interested in families of
frames along certain submanifolds. Such a family is usually called an
orthonormal
 moving frame. Throughout the text, we refer to it briefly as (moving) frame.
A Frenet frame is an example of a moving frame.

\begin{Exam}[Frenet frames along a curve] Let $\alpha=\alpha(t)$ be a curve in
$E^3$. The Frenet frame $\{x,e_1,e_2,e_3\}$ along the curve $\alpha$, as
described
in classical differential geometry, consists of the unit tangent vector
field $e_1$, the
unit normal vector field $e_2$ and the unit binormal vector field $e_3$. These vectors
satisfy the Frenet equations$$
\begin{cases}
\d x=\d s\cdot e_1&\\
\d e_1=\d s\cdot k(t)\cdot e_2&\\
\d e_2=\d s\cdot(-k(t)\cdot e_1+\tau(t)e_3)&\\
\d e_3=-\d s\cdot\tau(t)\cdot e_2&
\end{cases}
\eqno(1.1.1)$$

Here $\d s=s'(t)\d t$ represents the arc length differential, while
$k$ and
$\tau$ denote the curvature and torsion, respectively. Conversely, given
arbitrary
differential forms $\d s\ne0$, $k(t)\d t$, $\tau(t)\d t$, one can
reconstruct the
curve uniquely up to Euclidean motions.
\end{Exam}

For moving frames of surfaces, there exist differential forms generalizing $\d
s$, $k(t)\d t$, $\tau(t)\d t$, satisfying some integrability conditions, the
Gauss-Codazzi equations. Cartan showed that these equations can be derived from
the integrability conditions satisfied by the so-called Cartan forms
(see (1.1.10--11) below).

We will see that the space of all frames $\cal F$ forms a 6-dimensional
manifold.
This manifold can be identified with the group of Euclidean motions defined
below.

Consider the groups
\begin{align}
{\rm O}(3)&=\{A:E^3\to E^3\text{ linear; }\langle Ax,Ay\rangle=\langle
x,y\rangle,\ x,y\in E^3\}\tag{1.1.2}\\
{\rm SL}(3,{\Bbb R})&=\{A:E^3\to E^3\text{ linear; }\det A=1\}\tag{1.1.3}\\
{\rm SO}(3)&=\{A\in{\rm O}(3);\ \det A>0\}.\tag{1.1.4}
\end{align}

Note, ${\rm SO}(3)={\rm SL}(3,{\Bbb R})\cap {\rm O}(3)$.

We define the group of orientation-preserving rigid motions$$
G=\{w\mapsto x+Aw;\ x\in E^3,\ A\in{\rm SO}(3)\}.\eqno(1.1.5)$$
Note that the groups (1.1.2)--(1.1.5) are real Lie groups.

To identify G with $\cal F$, we fix a frame $F_0=\{0,\check{e}_1,\check{e}_2,
\check{e}_3\}$ in $\cal F$. Then if $F=\{x,e_1,e_2,e_3\}$ is an arbitrary
frame in
$\cal F$, the map$$
w\mapsto x+\sum_{i=1}^3\langle\check{e}_i,w\rangle e_i,\qquad w\in
E^3\eqno(1.1.6)$$
is an element of the group $G$.

Fixing $F_0$ means fixing an origin and an orthonormal basis. Expressing
the entries
of an arbitrary frame $F$ in terms of this basis, via (1.1.6),
realizes $F$ as a pair consisting of a translation vector and an
orientation-preserving matrix.

Conversely, given $g\in G$ we set$$
x=g(0)\quad\text{and}\quad e_i=g(\check{e}_i)-x.\eqno(1.1.7)$$

The resulting $F=\{x,e_1,e_2,e_3\}$ is a frame and it is easy to see that the
operations (1.1.6) and (1.1.7) are inverse to each other.

The bijection presented above gives an isomorphism between
$G$ and $\cal F$, and thus $\cal F$ is endowed with a manifold structure.

We consider the maps
\begin{align}
x^f:{\cal F}\to E^3,\quad &x^f(\{y,u_1,u_2,u_3\})=y,\tag{1.1.8}\\
e^f_j:{\cal F}\to E^3,\quad &e^f_j(\{y,u_1,u_2,u_3\})=u_j,\qquad
j=1,2,3\tag{1.1.9}
\end{align}

The differentials $\d x^f$, $\d e_1^f$, $\d e_2^f$ and $\d e_3^f$ estimated
at $F$ are
linear maps from the tangent space $T_F{\cal F}$ of $\cal F$ to $E^3$.
Therefore they
can be written as linear combinations relative to the basis
$e_1^f(F),e_2^f(F),e_3^f(F)$.

Thus, at a ``point" $F\in{\cal F}$, we define the scalar differential forms
$\omega_1,\omega_2,\omega_3$, $\omega_{ij}$, $1\leq i,j\leq3$ via$$
\d_Fx^f=\omega_1e_1^f(F)+\omega_2e_2^f(F)+\omega_3e_3^f(F)\eqno(1.1.10)$$
and$$
\d_Fe^f_j=\sum_{k=1}^3\omega_{jk}e_k^f(F).\eqno(1.1.11)$$
Since $\langle e_j,e_j\rangle=1$, we have $\langle e_j^f(F),e_j^f(F)\rangle=1$,
for every
$F\in\cal F$. Therefore,$$
\langle\d_Fe_j^f({\cal U}),e_j^f(F)\rangle=0,\qquad\text{for all }{\cal
U}\in T_F\cal
F.\eqno(1.1.12)$$
This last relation implies$$
\omega_{jj}=0, \qquad j=1,2,3\eqno(1.1.13)$$
Moreover, since $\langle e_j^f(F),e_k^f(F)\rangle=0$ for all $j\ne k$,
differentiation yields$$
\langle \d_Fe_j^f({\cal U}),e_k^f(F)\rangle+
\langle e_j^f(F),\d_Fe_k^f({\cal U})\rangle=0,\eqno(1.1.14)$$
that is$$
\omega_{jk}=-\omega_{kj}.\eqno(1.1.15)$$

Therefore, equations (1.1.10--11) are completely determined by the six 1-forms
$\omega_1,\omega_2,\omega_3,\omega_{12},\omega_{13},\omega_{23}$.

It is straightforward to verify the Cartan Structure Equations ([Ch, Te],
p.106):
\begin{align}
\d\omega_i&=\sum_{j=1}^3\omega_j\wedge\omega_{ji},\qquad i=1,2,3\tag{1.1.16}\\
\d\omega_{ij}&=\sum_{k=1}^3\omega_{ik}\wedge\omega_{kj},\qquad 1\leq
i,j\leq3.\tag{1.1.17}
\end{align}

Let now $M=(D,\psi)$ be an immersion of an open connected subset
$D\subset\Bbb R^2$
into $\Bbb R^3$, $\psi:D\to\Bbb R^3$. This describes a parametric surface,
admitting self-intersections.

All the frames $\{x,e_1,e_2,e_3\}$ with $x\in \psi(D)$ form the {\it zeroth
order
frame bundle} of $M$. The set of zeroth order frames will be denoted by
${\cal F}_0^{\cal M}$.

It is easy to see that the diffeomorphism ${\cal F}\cong G\cong\Bbb
R^3\times{\rm SO}(3)$, induced by fixing a frame $F_0$, yields
${\cal F}_0^{\cal M}\cong\psi(D)\times{\rm SO}(3)$.

Let now $F\in{\cal F}_0^{\cal M}$. Since we identified $\cal F$
with $G$, the group of orientation-preserving rigid motions of
$\Bbb R^3$, the frame $F$ in particular is identified with a pair
$F=(\psi(u,v),A)$, $(u,v)\in D$, $A\in{\rm SO}(3)$.

On ${\cal F}_0^{\cal M}$ we have natural vector fields:
$Y_1=\partial_uF$, $Y_2=\partial_vF$ and $Y_B$, where $B$ is any vector field
of $SO(3)$, pulled back to ${\cal F}_0^{\cal M}$. Then $\omega_j(Y_B)=0$.

Let $S^2$ denote the unit sphere in $\Bbb R^3$.
Let $\vec{N}:D\to S^2$ denote a unit normal vector field to $M$.
Then
$\langle\vec{N}(u,v),\partial_u\psi(u,v)\rangle=0$.
If $Y_1=\partial_uF$, $Y_2=\partial_vF$ denote the standard vector fields
along ${\cal F}_0^{\cal M}$ introduced above, then
\begin{align}
0&=\langle\vec{N}(u,v),\partial_u\psi(u,v)\rangle=\nonumber\\
&=\langle\vec{N},\d_Fx^f(Y_1)\rangle=\sum_j\langle\vec{N},e_j^f(F)\rangle\omega_
j(Y_1).
\tag{1.1.18}
\end{align}
Similarly, we obtain$$
0=\sum_j\langle\vec{N},e_j^f(F)\rangle\omega_j(Y_2).\eqno(1.1.19)$$
Since $\omega_j(Y_B)=0$ for all $Y_B$, restricting
$\omega_1,\omega_2,\omega_3$
to ${\cal F}_0^{\cal M}$, we obtain$$
\sum_j\langle\vec{N},e_j^f(F)\rangle\omega_j=0.\eqno(1.1.20)$$
Relation (1.1.20) represents the equation of the tangent plane to $M$ at $x$
relative to the frame $F=\{x,e_1,e_2,e_3\}$. The coefficients
$a_j:=\langle\vec{N},e_j^f(F)\rangle$ vary smoothly with the frame. Note
that if the frame is such that $e_1,e_2$ span the tangent plane of $M$ at
$x$, then
above linear relation (1.1.20) takes the form $\omega_3=0$.

For our goals, it is natural to consider moving frames for which $e_1$ and
$e_2$ are tangent to $M$.

\begin{Def}
Given an immersion $M=(D,\psi)$ as above, we define$$
{\cal F}_1^{\cal M}=\{(x,e_1,e_2,e_3)\in{\cal F}_0^{\cal M};\ e_1,e_2\in
T_xM\},\eqno(1.1.21)$$
where $T_xM$ denotes the tangent plane to $M$ at $x$.

${\cal F}_1$ is called the {\it first order frame bundle} of $M$.
\end{Def}

Along ${\cal F}_1$, $\omega_3$ vanishes, that is$$
\omega_3|_{T{\cal F}_1}=0.\eqno(1.1.22)$$
The above relation also implies$$
0=\d\omega_3=\omega_1\wedge\omega_{13}+\omega_2\wedge\omega_{23}\qquad\text{on
}T{\cal F}_1\times T{\cal F}_1.\eqno(1.1.23)$$
For ${\cal F}_1$, Cartan's structure equations (1.1.16--17) are written as
\begin{align}
d\omega_1=\omega_{12}\wedge\omega_2\tag{1.1.24a}\\
d\omega_2=\omega_1\wedge\omega_{12}\tag{1.1.24b}\\
d\omega_{12}=-\omega_{13}\wedge\omega_{23}\tag{1.1.24c}\\
d\omega_{13}=\omega_{12}\wedge\omega_{23}\tag{1.1.24d}\\
d\omega_{23}=\omega_{13}\wedge\omega_{12}\tag{1.1.24e}\\
\omega_3=0\tag{1.1.24f}\\
\omega_1\wedge\omega_{13}+\omega_2\wedge\omega_{23}=0\tag{1.1.24g}
\end{align}

The third equation above is also known as the {\it Gauss equation}, while the
fourth and fifth together are known as the {\it Codazzi equations}.

By Cartan's lemma [Ca, p.61], the last equation of (1.1.24g) implies:$$
\begin{cases}
\omega_{13}=h_{11}\omega_1+h_{12}\omega_2,&\\
\omega_{23}=h_{12}\omega_1+h_{22}\omega_2.&
\end{cases}\eqno(1.1.25)$$
for some functions $h_{ij}$ defined on $D$.

By (1.1.25) and the Gauss equation (1.1.24c), we obtain$$
\d\omega_{12}=-\omega_{13}\wedge\omega_{23}=-[\det(h_{ij})]\omega_1\wedge\omega_2,
\eqno(1.1.26)$$
where $K\buildrel{\rm def}\over=\det h_{ij}$ is called the {\it Gaussian
curvature} of $M$.

For the immersion $M=(D,\psi)$ and the submanifold ${\cal F}_1^{\cal
M}\subset{\cal
F}$, the
first fundamental form becomes$$
{\rm I}=\langle\d x^f,\d x^f\rangle=\omega_1^2+\omega_2^2,\eqno(1.1.27)$$
while the second fundamendal form is
\begin{align}
{\rm II}&=\langle-\d N,\d x^f\rangle=\langle-\d e_3^f,\d x^f\rangle\nonumber\\
&=\left\langle-\sum_{k=1}^3\omega_{3k}e_k^f,\sum_{i=1}^3\omega_ie_i^f\right\rangle\nonumber\\
&=-(\omega_{31}\omega_1+\omega_{32}\omega_2)=\omega_{13}\omega_1+\omega_{23}\omega_2.
\tag{1.1.28}
\end{align}
In formula (1.1.28), we chose the normal unit vector $N=e_3$.

Taking into account equations (1.1.25), we obtain$$
{\rm
II}=h_{11}\omega_1^2+2h_{12}\omega_1\omega_2+h_{22}\omega_2^2.\eqno(1.1.29)$$
The two-form $\omega_1\wedge\omega_2$ is an area element for the surface.
Therefore,
since $\psi: D\to\Bbb R^3$ is an immersion, it follows that$$
\omega_1\wedge\omega_2\ne0.\eqno(1.1.30)$$
As a consequence of formulas (1.1.25), the Gaussian curvature
$K=h_{11}h_{22}-h_{12}^2$ is given by$$
\omega_{13}\wedge\omega_{23}=(h_{11}\omega_1+h_{12}\omega_2)\wedge(h_{12}\omega_
1+
h_{22}\omega_2)=K(\omega_1\wedge\omega_2),\eqno(1.1.31)$$
while the mean curvature $H=(h_{11}+h_{22})/2$ is given by$$
\omega_1\wedge\omega_{23}-\omega_2\wedge\omega_{13}=h_{22}(\omega_1\wedge\omega_
2)
-h_{11}(\omega_2\wedge\omega_1)=2H(\omega_1\wedge\omega_2).\eqno(1.1.32)$$

All the formulas presented in this section were formulated by Elie Cartan
([Ca]). We followed the presentation of Cartan's structure equations for
$\Bbb R^3$
in [Ch, Te, eqs. (1.1)--(1.3)] and the one for surfaces in space forms of
constant
Gaussian curvature from [Te, eqs. (1.1)--(1.14)].

\section{Pseudospherical Surfaces and the Sine-Gordon Equation}

In this section, we begin our study of surfaces with constant negative
Gaussian curvature.
Among them, surfaces of Gaussian curvature $K=-1$, called {\it
pseudospherical surfaces},
are of particular interest to us. We show that all surfaces with constant
Gaussian curvature
are described by a sine-Gordon equation, and we write a corresponding Lax
system.

The following two parametrizations are of significant importance for this
class of surfaces. We will also specify the relationship between the
parametrizations.

\subsection{The Asymptotic Line Parametrization}

Let us consider an immersion $M=(D,\psi)$ with constant negative Gaussian
curvature. In the Euclidean space, every unit free vector represents a
{\it direction}.

For each point of $M$, there are two directions in which the second
fundamental form
vanishes, called {\it asymptotic directions} ([Ei, (46.3)]). An {\it
asymptotic line}
on the surface $M$ is a regular connected curve whose tangent unit vector is
an asymptotic direction at each point. Consequently, we have two families
of asymptotic lines, each tangent to an asymptotic direction everywhere.
An {\it asymptotic line parametrization} is a parametrization such that the
coordinate lines are asymptotic lines.

The given immersion $M=(D,\psi)$ can be locally reparametrized, such that
the coordinate lines are asymptotic lines. For an open and connected domain
$D$,
this reparametrization can be done globally. Therefore, for the rest of
this section
we will assume $\psi:D\to\Bbb R^3$ to be an asymptotic line parametrization
of the
surface $M$, where $D$ is a open connected domain in $\Bbb R^2$.

Let $\varphi$ represent the angle between the asymptotic lines,
measured counterclockwise from the vector field $\psi_x$ to the vector
field $\psi_y$.

We denote $A=|\psi_x|$, $B=|\psi_y|$.

Then the first fundamental form is ([Ei], [Bo2]):$$
{\rm I}=|\d\psi|^2=A^2(\d x)^2+2AB\cos\varphi\d x\d y+B^2(\d y)^2.$$
For every point, via a change of coordinates, we can reparametrize the
surface such that the asymptotic lines are parametrized in arc length.

Let us assume that $A$ and $B$ never vanish. An immersion $\psi$ with this
property is called {\it weakly regular}. A weakly regular surface can be always
reparametrized such that both asymptotic lines are in arc length
($A=B=1$).

In this context, let $N:D\to S^2$,
$N=\frac{{\psi_x}\times{\psi_y}}{\|{\psi_x}\times{\psi_y}\|}$ define the
Gauss map
of the immersion $\psi$.
Remark that the unit vector field $N$ is orthogonal to $\psi_x$, $\psi_y$,
$\psi_{xx}$,
$\psi_{yy}$.

\begin{Def}
A parametrization for which $A=B=1$ is called a Chebyshev net ([Spi]).
\end{Def}

Unless stated otherwise, we will assume for the rest of this work
that the immersion $\psi$ corresponds to a Chebyshev net of angle
(between asymptotic lines) $\varphi(x,y)\in(0,\pi)$. In this case,
the metric becomes:$$ {\rm I}=|d\psi|^2=(\d x)^2+2\cos\varphi\d
x\d y+(\d y)^2.\eqno(2.1.1)$$ [McL] presents a way of constructing
a Chebyshev net physically, by ``a piece of nonstrech fabric that
is loosely woven, so that the angle between the threads can
change. Then drape it over the surface so that the warp and weft
of the fabric become coordinate lines on the surface''. Since the
threads cannot stretch, $A=B=1$, but the angle $\varphi(x,y)$
changes. The second fundamental form in asymptotic parametrization
is written as$$ {\rm II}=2 A B \sqrt{-K}\sin\varphi\,\d x\,\d y.$$
For a Chebyshev net, it clearly becomes
 $$
{\rm II}=2 \sqrt{-K}\sin\varphi\,\d x\,\d y,\eqno(2.1.2)$$
where $K$ represents the (constant, negative) Gaussian curvature, $$
K=\det{\rm II}/\det{\rm I}.$$
Let us now focus on the case of the pseudospherical surfaces, that is
surfaces of
Gaussian curvature $K=-1$. It is straightforward to calculate the principal
curvatures $k_1$ and $k_2$ of the immersion. $k_1$ and $k_2$ represent the
eigenvalues
of the matrix$$
{\rm II}\cdot{\rm I}^{-1}=\begin{pmatrix}
-\cot\varphi&{\rm csc}\,\varphi\\
{\rm csc}\,\varphi&-\cot\varphi
\end{pmatrix},\eqno(2.1.3)$$
that is, the roots of the characteristic equation$$
\lambda^2+2\cdot\cot\varphi\cdot\lambda-1=0,$$
i.e.,$$
k_1=\tan\frac{\varphi}{2}\quad\text{and}\quad
k_2=-\cot\frac{\varphi}{2}.\eqno(2.1.4)$$
The angle between the asymptotic lines can be written as
$\varphi(x,y)=2\arctan
k_1$.

Let $e_1$ and $e_2$ be the principal directions on $M$ corresponding to
$k_1$ and
$k_2$ respectively, that is the eigenvectors of the matrix ${\rm II}\cdot{\rm
I}^{-1}$ at each point of $M$. Then the relation between the asymptotic
directions
on $M$ and the principal directions on $M$ is given by
\begin{align}
\partial_x&=\cos\frac{\varphi}{2}e_1-\sin\frac{\varphi}{2}e_2,\nonumber\\
\partial_y&=\cos\frac{\varphi}{2}e_1+\sin\frac{\varphi}{2}e_2.\tag{2.1.5}
\end{align}

\subsection{The Curvature Line Parametrization; Sine-Gordon Equation}

Another useful parametrization for a pseudospherical immersion
$M=(D,\psi)$  is the one by lines of curvature, i.e., the coordinates $u_i$
in which
both the first fundamental form I and the second fundamental form II are
diagonalized as

\begin{align}
{\rm I}&=(a_1)^2(\d u_1)^2+(a_2)^2(\d u_2)^2\tag{2.2.1a}\\
{\rm II}&=b_1\cdot(a_1)^2(\d u_1)^2+b_2\cdot(a_2)^2(\d u_2)^2.\tag{2.2.1b}
\end{align}
In general, such a parametrization exists only in the neighborhood of a
non-umbilical
point. Since the Gaussian curvature is negative, there are no umbilics on $M$.

In particular, on a weakly regular pseudospherical surface we can find a
curvature
line parametrization around every point.

More specifically, we set
 $$
u_1=x+y, \qquad u_2=x-y,$$

where $(x,y)$ are the Chebyshev net coordinates from Section 2.1. (i.e.
$A=B=1$).

Then formulas (2.1.1) and (2.1.2), for $K=-1$, become:
\begin{align}
{\rm I}&=\cos^2\frac{\varphi}{2}\cdot(\d u_1)^2+\sin^2\frac{\varphi}{2}\cdot(\d
u_2)^2\tag{2.2.2a}\\
{\rm II}&=\sin\frac{\varphi}{2}\cos\frac{\varphi}{2}((\d
u_1)^2-(\d u_2)^2)\tag{2.2.2b}
\end{align}
respectively.

Comparing with (2.2.1) above, we obtain:
\begin{align}
a_1&=\cos\frac{\varphi}{2},\tag{2.2.3a}\\
a_2&=\sin\frac{\varphi}{2},\tag{2.2.3b}\\
b_1&=k_1=\tan\frac{\varphi}{2},\tag{2.2.3c}\\
b_2&=k_2=-\cot\frac{\varphi}{2},\tag{2.2.3d}
\end{align}
where $\varphi(x,y)$ is the angle between the asymptotic directions and
$k_1$, $k_2$
represent the principal curvatures.

Note that (2.2.3 a-d) correspond to a choice of $a_1, a_2, b_1, b_2$ made
without
loss of generality ([Te], 2.7).

We also note that in asymptotic line parametrization, the principal vectors
given by (2.1.5) are generally not orthogonal, so the context is different
than the one
of orthonormal frames (Section 1). However, in curvature line coordinates,
the principal vectors $e_1$ and $e_2$ are orthogonal, and that enables us
to use
the moving frame context from Section 1.

Comparing formulas (2.2.2) to the formulas (1.1.27) and (1.1.29), we deduce:
\begin{align}
\omega_1&=a_1\d u_1=\cos\frac{\varphi}{2}\d u_1,\tag{2.2.4a}\\
\omega_2&=a_2\d u_2=\sin\frac{\varphi}{2}\d u_2,\tag{2.2.4b}\\
h_{11}&=k_1,\quad h_{12}=0,\quad h_{22}=k_2.\tag{2.2.4c}
\end{align}
Then (1.1.25) together with (2.2.4c) yield
\begin{align}
\omega_{13}&=k_1a_1\d u_1,\tag{2.2.5a}\\
\omega_{23}&=k_2a_2\d u_2.\tag{2.2.5b}
\end{align}

We also aim at finding an expression for $\omega_{12}$: from equations
 (2.2.4a) and (2.2.4b), we find:
\begin{align}
\d\omega_1&=\frac{\partial a_1}{\partial u_2}\d u_2\wedge\d
u_1=-\frac{1}{a_2}\cdot
\frac{\partial a_1}{\partial u_2}\d u_1\wedge\omega_2,\tag{2.2.6a}\\
\d\omega_2&=\frac{\partial a_2}{\partial u_1}\d u_1\wedge\d
u_2=\frac{1}{a_1}\cdot
\frac{\partial a_2}{\partial u_1}\omega_1\wedge\d u_2.\tag{2.2.6b}
\end{align}
Comparing equation (2.2.6) to the first two structure equations, (1.1.24a)
and (1.1.24b), we obtain$$
\omega_{12}=\frac{1}{a_1}\frac{\partial a_2}{\partial u_1}\d u_2-\frac{1}{a_2}
\frac{\partial a_1}{\partial u_2}\d u_1.\eqno(2.2.7)$$
As a consequence of (2.2.5a,b) and (2.2.7), we deduce
\begin{align}
&\omega_{12}\wedge\omega_{23}=-k_2\frac{\partial a_1}{\partial u_2}\d
u_1\wedge\d
u_2\tag{2.2.8a}\\
&\d\omega_{13}=\d(k_1\omega_1)=\left(-k_1\frac{\partial a_1}{\partial
u_2}-a_1\frac{\partial k_1}{\partial u_2}\right)\d u_1\wedge\d u_2.\tag{2.2.8b}
\end{align}
Therefore, the first Codazzi equation, (1.1.24d), has the form
$$
(k_2-k_1)\frac{\partial a_1}{\partial u_2}=a_1\frac{\partial k_1}{\partial
u_2},$$
which can be rewritten as$$
\frac{1}{k_2-k_1}\frac{\partial k_1}{\partial u_2}=\frac{\partial(\log
a_1)}{\partial u_2}.\eqno(2.2.9{\rm a})$$
Similarly, the second Codazzi equation, (1.1.24e), becomes$$
\frac{1}{k_1-k_2}\frac{\partial k_2}{\partial u_1}=\frac{\partial(\log
a_2)}{\partial u_1}.\eqno(2.2.9{\rm b})$$

Recall now that $\psi$ is a Chebyshev net parametrization:
$A=|\psi_x| = 1$ and $B=|\psi_y| = 1$. In general (see, e.g., [Bo2], p. 114),
the Codazzi equation can be written as$$
A_y=B_x=0.\eqno(2.2.10)$$
So the Codazzi equations become trivial for a Chebyshev net.

Let us focus now on the Gauss equation (1.1.24c):$$
\d\omega_{12}=-\omega_{13}\wedge\omega_{23}.$$

Substituting the expressions for $a_1$ and $a_2$ from (2.2.3) into (2.2.7), we
obtain the following expression for the connection form $\omega_{12}$:$$
\omega_{12}=\frac{1}{2}\left(\frac{\partial\varphi}{\partial u_1}\d
u_2+\frac{\partial\varphi}{\partial u_2}\d u_1\right).\eqno(2.2.11)$$
Therefore,$$
\d\omega_{12}=\frac{1}{2}\left(\frac{\partial^2\varphi}{\partial(u_1)^2}-
\frac{\partial^2\varphi}{\partial(u_2)^2}\right)\d u_1\wedge\d
u_2.\eqno(2.2.12)$$

Further, substituting the expressions (2.2.3) for $a_1,a_2,k_1,k_2$ into
(2.2.5), the Gauss equation (1.1.24c) can be written in curvature coordinates
as
$$
\frac{\partial^2\varphi}{\partial(u_1)^2}-
\frac{\partial^2\varphi}{\partial(u_2)^2}=\sin\varphi.\eqno(2.2.13)$$
Via $u_1=x+y, u_2=x-y$, (2.2.13) becomes, in asymptotic line parametrization,$$
\varphi_{xy}=\sin\varphi,\eqno(2.2.14)$$
Note that (2.2.13) and (2.2.14) are two different forms of the sine-Gordon
equation.

Conversely, by the existence and uniqueness theorem of surface theory,
given $\varphi$, a solution to (2.2.14), there exists an immersion
$M=(D,\psi)$, in asymptotic line coordinates, whose angle between asymptotic
directions is $\varphi$.

Summarizing the discussion above, we can state now the following result,
due to Enneper (1845):

\begin{Th} {\rm([Ch], p. 441, and [Bo2], p. 115)}\quad Up to rigid motion,
there is a one-to-one correspondence between solutions $\varphi$ to the
sine-Gordon equation (2.2.14) with $0<\varphi<\pi$ and the weakly regular
pseudospherical surfaces in Chebyshev net parametrization immersed in $E^3$.
\end{Th}

\noindent{\bf Note.} This one-to-one correspondence between solutions
$\varphi$ to the sine-Gordon equation (2.2.14) and pseudospherical surfaces,
whose first and second fundamental forms are given by (2.2.2), is the
particular case $K<\bar{K}=0$ of the following general theorem ([Te], Cor.2.7):

\begin{Th} Let $M^2(K)$ be a surface with constant Gaussian curvature $K$,
contained in a Riemannian $3$-dimensional space form $\bar{M}^3(\bar{K})$
with constant curvature $\bar{K}$ such that $K\ne\bar{K}$. If $K>\bar{K}$,
assume that $M$ has no umbilic points. Then there exist local coordinates
$x_1,x_2$ and a real-valued function $\psi(x_1,x_2)$ which satisfies the
differential equation
\begin{align}
\psi_{x_1x_1}-\psi_{x_2x_2}&=-K\sin\psi\qquad\text{if }K<\bar{K},\tag{$\ast$}\\
\psi_{x_1x_1}+\psi_{x_2x_2}&=-K\sinh\psi\qquad\text{if
}K>\bar{K}.\tag{$\ast\ast$}
\end{align}
Conversely, suppose $\psi$ is a solution of $(\ast)$ (resp. $(\ast\ast)$). Then
there exists a surface of constant Gaussian curvature $K$ in a space form
$\bar{M}^3(\bar{K})$, which is unique up to rigid motion of $\bar{M}^3$,
whose first
and second fundamental forms are given respectively by
\begin{align*}
{\rm I}&=
\begin{cases}
\cos^2\dis\frac{\psi}{2}\d x^2_1+\sin^2\dis\frac{\psi}{2}\d x^2_2&\mbox{if
}K<\bar{K},\\
\cosh^2\dis\frac{\psi}{2}\d x^2_1+\sinh^2\dis\frac{\psi}{2}\d x^2_2&\mbox{if
}K>\bar{K},
\end{cases}\\
{\rm II}&=
\begin{cases}
\sqrt{|K-\bar{K}|}\sin\dis\frac{\psi}{2}\cos\dis\frac{\psi}{2}(\d x_1^2-\d
x_2^2)&\mbox{if }K<\bar{K},\\
\sqrt{|K-\bar{K}|}\sinh\dis\frac{\psi}{2}\cosh\dis\frac{\psi}{2}(\d x_1^2+\d
x_2^2)&\mbox{if }K>\bar{K}.
\end{cases}\\
\end{align*}
\end{Th}

\subsection{Moving Frame of a Pseudospherical Surface. The Lax System}

Let $D$ be a simply connected domain in $\Bbb R^2$ and $\psi:D\to\Bbb R^3$ an
immersion corresponding to a pseudospherical surface $M=(D, \psi)$.
Let $k_1$, $k_2$ be the principal curvatures, given by formulas (2.2.3c,d)
and $e_1$, $e_2$ corresponding principal directions
on $M$. Let
$F=\{x,e_1,e_2,e_3\}\in{\cal F}_1^{\cal M}$ be a fixed moving frame. Clearly,
$e_3$ represents a chosen normal direction $N$, along $M$. Let us focus now
on the Frenet equations of the frame. We shall omit the component $x\in E^3$
and will identify $F=\begin{pmatrix}e_1\\e_2\\e_3\end{pmatrix}$ for the rest
of this section.
By (1.1.11), we have the following Frenet system on $M$:$$
\d F=\begin{pmatrix}
0&\omega_{12}&\omega_{13}\\
-\omega_{12}&0&\omega_{23}\\
-\omega_{13}&-\omega_{23}&0
\end{pmatrix}F.\eqno(2.3.1)$$
The 1-forms $\omega_{12}$, $\omega_{13}$ and $\omega_{23}$, as a consequence of
formulas (2.2.4--5) and (2.2.11), can be written as
\begin{align}
\omega_{12}&=\frac{1}{2}(\varphi_{u_1}\d u_2+\varphi_{u_2}\d u_1)=\frac{1}{2}
(\varphi_x\d x-\varphi_y\d y)\tag{2.3.2a}\\
\omega_{13}&=k_1\omega_1=\sin\frac{\varphi}{2}\cdot\d
u_1=\sin\frac{\varphi}{2}\cdot
(\d x+\d y)\tag{2.3.2b}\\
\omega_{23}&=k_2\omega_2=-\cos\frac{\varphi}{2}\cdot\d
u_2=-\cos\frac{\varphi}{2}\cdot (\d x-\d y)\tag{2.3.2c}
\end{align}

Let us now consider the moving frame
$\tilde{F}_{\theta}\in{\cal F}_1^{\cal M}$, that is obtained from $F$
via a rotation of angle $\theta(x,y)$ in the tangent plane,
around $N$,
namely$$
\tilde{F}_{\theta}=
\begin{pmatrix}\tilde{e}_1\\\tilde{e}_2\\N\end{pmatrix},\eqno(2.3.3)$$
where$$
\begin{pmatrix}\tilde{e}_1\\\tilde{e}_2\end{pmatrix}=\begin{pmatrix}
\cos\theta&\sin\theta\\
-\sin\theta&\cos\theta\end{pmatrix}\begin{pmatrix}e_1\\e_2\end{pmatrix}.$$
In particular for $\theta=\varphi/2$, where $\varphi(x,y)$ is the angle
between the
asymptotic directions, the resulting frame is denoted $\tilde{F}$ and is
called
the {\it normalized
frame} associated with the moving frame $F$ (see [Wu1], p.18).
Unless stated otherwise, we will denote by $F$ the usual coordinate frame,
and by $\tilde{F}$ the rotated frame as stated above. A simple calculation
leads us to
the system of Frenet equations for $\tilde{F}$:$$
\d\tilde{F}=\begin{pmatrix}0&\tilde{\omega}_{12}&\tilde{\omega}_{13}\\
-\tilde{\omega}_{12}&0&\tilde{\omega}_{23}\\
-\tilde{\omega}_{13}&-\tilde{\omega}_{23}&0\end{pmatrix}\tilde{F},\eqno(2.3.4)$$
where$$
\tilde{\omega}_{12}=\d\theta+\omega_{12},\eqno(2.3.5a)$$$$
\tilde{\omega}_{13}=k_1\cos\theta\,\omega_1+k_2\sin\theta\,\omega_2,\eqno(2.3.5b
)$$$$
\tilde{\omega}_{23}=-k_1\sin\theta\,\omega_1+k_2\cos\theta\,\omega_2.\eqno(2.3.5
c)$$

In particular for the normalized frame $\tilde{F}$, $\theta=\varphi/2$
implies:$$
\d\theta=\frac{1}{2}(\varphi_x\d x+\varphi_y\d y)\eqno(2.3.6a)$$$$
\omega_{12}=\frac{1}{2}(\varphi_x\d x-\varphi_y\d y)\eqno(2.3.6b)$$$$
\omega_1=\cos\frac{\varphi}{2}(\d x+\d y)\eqno(2.3.6c)$$$$
\omega_2=\sin\frac{\varphi}{2}(\d x-\d y)\eqno(2.3.6d)$$$$
\tilde{\omega}_{12}=\varphi_x\d x\eqno(2.3.6e)$$$$
\tilde{\omega}_{13}=\sin\frac{\varphi}{2}\omega_1-\cos\frac{\varphi}{2}\omega_2
=
\frac{\sin\varphi}{2}(\d u_1-\d u_2)=\sin\varphi\cdot\d y\eqno(2.3.6f)$$$$
\tilde{\omega}_{23}=-\sin^2\frac{\varphi}{2}\d
u_1-\cos^2\frac{\varphi}{2}\d u_2=
-\d x+\cos\varphi\cdot\d y.\eqno(2.3.6g)$$
As a consequence of (2.3.6), the Frenet system (2.3.4) is equivalent to the
following differential system (also called {\it Lax system}):
$$
\begin{aligned}
\partial_x\tilde{F}&=\begin{pmatrix}0&\varphi_x&0\\
-\varphi_x&0&-1\\
0&1&0\end{pmatrix}\tilde{F}=\tilde{\cal A}\tilde{F},\nonumber\\
\partial_y\tilde{F}&=\begin{pmatrix}0&0&\sin\varphi\\
0&0&\cos\varphi\\
-\sin\varphi&-\cos\varphi&0\end{pmatrix}\tilde{F}=\tilde{\cal B}\tilde{F}.
\end{aligned}\eqno(2.3.7)$$

Note that $\tilde{\cal A}$ and $\tilde{\cal B}$ are skew-symmetric matrices.

The compatibility condition for the system (2.3.7) (i.e.
$\tilde{F}_{xy}=\tilde{F}_{yx}$) is
$$
\tilde{\cal A}_y-\tilde{\cal B}_x-[\tilde{\cal A},\tilde{\cal
B}]=0.\eqno(2.3.8)$$
This is equivalent to the Gauss equation, which for pseudospherical
surfaces in
a Chebyshev parametrization is the sine-Gordon equation (2.2.14). If, for a
pseudospherical surface, we use any asymptotic line parametrization $\psi$,
but not
necessarily a Chebyshev net, the Gauss equation takes the more general form
([Bo2], p. 114):$$
\varphi_{xy}=AB\sin\varphi,\eqno(2.3.9)$$
where $A=|\psi_x|$, $B=|\psi_y|$.

It is interesting to remark that this equation remains invariant with
respect to the
transformation$$
A\mapsto\lambda A,\quad B\mapsto B/\lambda,\qquad\lambda\in\Bbb
R_+,\eqno(2.3.10)$$
which plays an essential role in the theory of pseudospherical surfaces.
The transformation (2.3.10) appears in literature as Lie's transformation or
Lorentz transformation in plane. To reconcile the two names, it is sometimes
called Lie-Lorentz transformation.

The following obvious result is due to Lie (around the year 1870) and is of
crucial
importance in our context ([Bo2], p. 114):

\begin{Th} Every surface with constant negative Gauss curvature has a
one-parameter
family of deformations preserving {\rm the second fundamental form}$$
{\rm II}=2AB\sqrt{-K}\sin\varphi\d x\d y,\eqno(2.3.11)$$
{\rm the Gaussian curvature} $K$ and {\rm the angle} $\varphi$ between the
asymptotic lines. The deformation is generated by the transformation
$(2.3.10)$ above.
\end{Th}

The family of immersions mentioned above is called {\it associated family of
surfaces}. It will be denoted as $\psi^{\lambda }:D\to \Bbb R^{3}$.
Note that all the immersions are defined on the same domain $D$.

\begin{Rem}The Lie-Lorentz transformation (2.3.10) can be naturally induced
by replacing
$x$ with $\lambda^{-1}x$ and $y$ by $\lambda y$, $\lambda>0$, and then$$
\begin{aligned}
\partial_x&=\lambda\left(\cos\frac{\varphi}{2}\cdot
e_1-\sin\frac{\varphi}{2}\cdot
e_2\right),\nonumber\\
\partial_y&=\frac{1}{\lambda}\left(\cos\frac{\varphi}{2}\cdot
e_1+\sin\frac{\varphi}{2}\cdot e_2\right).
\end{aligned}\eqno(2.3.12)$$
We note here that the Lie-Lorentz transformation defined above on
$M=(D,{\psi})$ is
equivalent to a Lorentz transformation on a Lorentzian 2-manifold, $(D, II)$.

Also note that if $\varphi(x,y)$ denotes the angle of a certain
pseudospherical surface $M$ in Chebyshev net coordinates $x,y$,
then by Lie-Lorentz transformation we create a new pseudospherical surface
$M^*$, in the same associated family with the first one. The coordinates
$x^*=\lambda^{-1} x$ and $y^*=\lambda y$ are also asymptotic, and the angle
between asymptotic lines on the new surface is given by the same function
as before, but this time in variables $x^*$ and $y^*$. Thought of as a
function of the {\it old} coordinates $x,y$, the angle $\varphi (x^*, y^*)$,
corresponding to the new surface $M^*$, depends on $\lambda$.
See also the examples in Section 8, (8.1.3) and (8.2.1).

\end{Rem}

As a consequence of the coordinate change described above via the parameter
$\lambda$, starting from a Chebyshev parametrization $\psi$, we see that
$|\psi_x|=1$ becomes $|\psi_x|=\lambda$, while $|\psi_y|=1$ becomes
$|\psi_y|=\lambda^{-1}$.
While via this transformation the sine-Gordon equation remains unmodified, the
corresponding differential Lax
system (2.3.7) depends on $\lambda$. In particular, we obtain an {\it
extended frame}
$F=F(x,y,\lambda)=F(\lambda^{-1} x, \lambda y)$. For the normalized frame
$\tilde{F}$, we obtain the {\it extended normalized frame}
$\tilde{F}(x,y,\lambda)$.

\begin{Cor} The extended normalized frame $\tilde{F}(x,y,\lambda)$
satisfies the following Lax differential system:
\begin{align}
\partial_x\tilde{F}&=\begin{pmatrix}
0&\varphi_x&0\\
-\varphi_x&0&-\lambda\\
0&\lambda&0\end{pmatrix}\tilde{F},\nonumber\\
\partial_y\tilde{F}&=\frac{1}{\lambda}\begin{pmatrix}
0&0&\sin\varphi\\
0&0&\cos\varphi\\
-\sin\varphi&-\cos\varphi&0\end{pmatrix}\tilde{F}.\tag{2.3.13}
\end{align}
\end{Cor}

This type of linear system is essential for the inverse scattering method
in soliton theory. Equation (2.3.13) represents the scattering system of
the sine-Gordon
equation introduced by Lund (see [Lu]).

\begin{Rem} The frame $F$ represents the $3\times3$ matrix
$\begin{pmatrix}e_1\\e_2\\e_3\end{pmatrix}$ of rows $e_1$, $e_2$, and $e_3$,
respectively. In the spirit of [Wu2] and [DoHa], instead of the classical
frame $F$, it is more convenient to work with ${\cal U} := \tilde F^T$,
the transposed of the extended normalized frame
$\tilde{F}(x,y,\lambda)$. This is especially convenient in view of
formulas (2.3.15) below. Unless stated otherwise, the term of normalized
coordinate
frame will refer to ${\cal U}$ above, for the rest of this text.

Consequently, formulas (2.3.7) can be rewritten as$$
\partial_x{\cal U}={\cal U}\cdot {\cal A}^T,\nonumber\\
\partial_y{\cal U}={\cal U}\cdot {\cal B}^T,\eqno(2.3.14)$$
where we denoted by ${\cal A}$ and ${\cal B}$, respectively, the transpose of
$\tilde{\cal A}$ and $\tilde{\cal B}$ from (2.3.7).

That is, equations (2.3.13) above can be rewritten as:

\end{Rem}

\begin{Cor} The extended normalized frame $\cal U^\lambda$
satisfies the
following Lax differential system
\begin{align}
\partial_x\cal U^\lambda&={\cal U}^\lambda\cdot\begin{pmatrix}
0&-\varphi_x&0\\
\varphi_x&0&\lambda\\
0&-\lambda&0\end{pmatrix},\nonumber\\
\partial_y\cal U^\lambda&={\cal U}^\lambda\frac{1}{\lambda}\begin{pmatrix}
0&0&-\sin\varphi\\
0&0&-\cos\varphi\\
\sin\varphi&\cos\varphi&0\end{pmatrix}.\tag{2.3.15}
\end{align}
\end{Cor}

The Lax system will be written in this form for the rest of this work.
It plays a crucial role in the study of pseudospherical surfaces.

\section{Associated Families of Pseudospherical Surfaces via Spectral
Parameter $\lambda$}

In this section we study in detail the effects of introducing the real
positive
parameter $\lambda$. We obtain in this way a $\lambda$-transformation of the
Cartan
forms (respectively an extended Maurer-Cartan form $\omega^\lambda$)
corresponding to the associated family of pseudospherical surfaces
(respectively the extended normalized frame ${\cal U}^\lambda$).

\subsection{The $\lambda$-Transformation on the 1-Forms $\omega_i$ and
$\omega_{ij}$}

Let us study the effect that the transformation (2.3.10) has on the
1-forms
$\omega_1,\omega_2,\omega_{12},\omega_{13},\omega_{23}$. Replacing $x$ by
$x^*:=\lambda^{-1}x$ and $y$ by $y^*:=\lambda y$ in the system
(2.3.2), and taking into account the invariance of $\varphi$ under this
deformation (Thm. 2.3.1), we obtain the ``extended'' forms:
\begin{align}
\omega_1^\lambda&=\cos\frac{\varphi}{2}(\d x^*+\d
y^*)=\cos\frac{\varphi}{2}(\lambda^{-1}\d x+\lambda\d y)\tag{3.1.1a}\\
\omega_2^\lambda&=\sin\frac{\varphi}{2}(\d x^*-\d
y^*)=\sin\frac{\varphi}{2}(\lambda^{-1}\d x-\lambda\d y)\tag{3.1.1b}\\
\omega_{12}^\lambda&=\frac{1}{2}(\varphi_{x^*}\d x^*-\varphi_{y^*}\d
y^*)=\frac{1}{2}(\varphi_x\d x-\varphi_y\d y)\tag{3.1.1c}\\
\omega_{13}^\lambda&=\sin\frac{\varphi}{2}(\d x^*+\d
y^*)=\sin\frac{\varphi}{2}(\lambda^{-1}\d x+\lambda\d y)\tag{3.1.1d}\\
\omega_{23}^\lambda&=-\cos\frac{\varphi}{2}(\d x^*-\d
y^*)=-\cos\frac{\varphi}{2}(\lambda^{-1}\d x-\lambda\d y).\tag{3.1.1e}
\end{align}

The system above can be rewritten as
\begin{align}
\omega_1^\lambda&=\frac{1}{2}(\lambda+\lambda^{-1})\omega_1+
\frac{1}{2}(\lambda-\lambda^{-1})\omega_{23},\tag{3.1.2a}\\
\omega_2^\lambda&=\frac{1}{2}(\lambda+\lambda^{-1})\omega_2-
\frac{1}{2}(\lambda-\lambda^{-1})\omega_{13},\tag{3.1.2b}\\
\omega_{12}^\lambda&=\omega_{12},\tag{3.1.2c}\\
\omega_{13}^\lambda&=-\frac{1}{2}(\lambda-\lambda^{-1})\omega_2+
\frac{1}{2}(\lambda+\lambda^{-1})\omega_{13},\tag{3.1.2d}\\
\omega_{23}^\lambda&=\frac{1}{2}(\lambda-\lambda^{-1})\omega_1+
\frac{1}{2}(\lambda+\lambda^{-1})\omega_{23},\tag{3.1.2e}
\end{align}
where $\lambda>0$. Note that $\lambda$ occurs rationally, with simple poles
at $\lambda = 0$ and at infinity. This will be essential below.

Cartan's structure equations for ${\cal F}_1^{\cal M}$, where $\omega_3$ is
identically zero, given by (1.1.24 a-f), together with equation (1.1.31)
for $K=-1$, form the set of equations below,
called {\it conditions} ($K$):
\begin{align}
&\d\omega_1=\omega_{12}\wedge\omega_2,\tag{3.1.3.a}\\
&\d\omega_2=\omega_1\wedge\omega_{12},\tag{3.1.3.b}\\
&\d\omega_{12}=-\omega_{13}\wedge\omega_{23},\tag{3.1.3.c}\\
&\d\omega_{13}=\omega_{12}\wedge\omega_{23},\tag{3.1.3.d}\\
&\d\omega_{23}=\omega_{13}\wedge\omega_{12},\tag{3.1.3.e}\\
&\omega_1\wedge\omega_{13}+\omega_2\wedge\omega_{23}=0,\tag{3.1.3.f}\\
&\omega_1\wedge\omega_2+\omega_{13}\wedge\omega_{23}=0.\tag{3.1.3.g}
\end{align}

Let $\omega_1,\omega_2,\omega_{12},\omega_{13},\omega_{23}$ be differential
forms
defined by (1.1.10) and let
$\omega_1^\lambda$, $\omega_2^\lambda,\omega_{12}^\lambda,\omega_{13}^\lambda,
\omega_{23}^\lambda$ be given by (3.1.2). Then

\begin{Th} The forms $\omega_1,\omega_2,\omega_{12},\omega_{13},\omega_{23}$ satisfy the
conditions (K) if and only if
$\omega_1^\lambda$, $\omega_2^\lambda,\omega_{12}^\lambda,\omega_{13}^\lambda,
\omega_{23}^\lambda$ satisfy the conditions (K).

For every pseudospherical surface $M=(D,\psi)$, there exists
a family $M_\lambda=(D,\psi_\lambda)$, $\lambda>0$, of pseudospherical surfaces
associated with
$\omega_1^\lambda,\omega_2^\lambda,\omega_{12}^\lambda$, $\omega_{13}^\lambda,
\omega_{23}^\lambda$ preserving the angle $\varphi$ between the asymptotic
lines and also preserving the second fundamental form.
\end{Th}

\noindent{\it Proof.} By Theorem 2.3.1, we know that the
$\lambda$-transformation (2.3.10) preserves the angle $\varphi$ and the
second fundamental form.
This means that the forms $\omega_i$, $\omega_{ij}$, and $\omega_i^\lambda$,
$\omega_{ij}^\lambda$, $\lambda>0$, respectively, satisfy the same Gauss
equation (3.1.3.c).
The Gauss equation is equivalent with $\varphi_{xy}=\sin\varphi$, and so
the angle $\varphi$
is preserved for the family $M_\lambda$. We remark that the Codazzi
equation is
trivially satisfied for $M_\lambda$, since for the whole associated family
$\psi^\lambda$,
$A=|\psi^\lambda_x|=\lambda$,
$B=|\psi^\lambda_y|=1/\lambda$,
$\lambda>0$, and the Codazzi equations are $A_y=B_x=0$.

In order to finish the proof of the theorem, it is enough to show that if
the Gauss
and Codazzi equations are satisfied for every real positive $\lambda$, then
the rest of
conditions $(K)$ are also satisfied for every real positive $\lambda$. This
is stated
in the following:

\begin{Lem} If $\omega_i^\lambda$ and $\omega_{ij}^\lambda$ are given by the
equations $(3.1.2)$, and if the following conditions are satisfied for all
$\lambda>0$:
\begin{align}
\d\omega_{12}^\lambda&=-\omega_{13}^\lambda\wedge\omega_{23}^\lambda,\tag{3.1.4.i}\\
\d\omega_{13}^\lambda&=\omega_{12}^\lambda\wedge\omega_{23}^\lambda,\tag{3.1.4.ii}\\
\d\omega_{23}^\lambda&=\omega_{13}^\lambda\wedge\omega_{12}^\lambda,\tag{3.1.4.iii}
\end{align}
then all the conditions (K) are satisfied for $\omega_i^\lambda$,
$\omega_{ij}^\lambda$.
\end{Lem}

\noindent{\it Proof.}

Assume that (3.1.4.i-iii) are satisfied. Then, by (3.1.2), after a few
simplifications, we obtain
\begin{align*}
\d\omega_{12}^\lambda+\omega_{13}^\lambda\wedge\omega_{23}^\lambda&={\lambda^2-
\lambda^{-2}\over4}(-\omega_1\wedge\omega_{13}-\omega_2\wedge\omega_{23})\\
&\quad+
{\lambda^2+\lambda^{-2}\over4}(\omega_1\wedge\omega_2+\omega_{13}\wedge\omega_{2
3})\\
&\quad+\frac{1}{2}(-\omega_1\wedge\omega_2-\omega_{13}\wedge\omega_{23})=0,\tag{
3.1.5.i}\\
\d\omega_{13}^\lambda-\omega_{12}^\lambda\wedge\omega_{23}^\lambda&={\lambda-
\lambda^{-1}\over2}(\d\omega_2-\omega_1\wedge\omega_{12})\\
&\quad+
{\lambda+\lambda^{-1}\over2}(\d\omega_{13}-\omega_{12}\wedge\omega_{23})=0,\tag{
3.1.5.ii}\\
\d\omega_{23}^\lambda-\omega_{13}^\lambda\wedge\omega_{12}^\lambda&=-{\lambda-
\lambda^{-1}\over2}(\d\omega_1-\omega_{12}\wedge\omega_2)\\
&\quad+
{\lambda+\lambda^{-1}\over2}(\d\omega_{23}-\omega_{13}\wedge\omega_{12})=0.\tag{
3.1.5.iii}
\end{align*}
Comparing the coefficients of the corresponding $\lambda^2$ and $\lambda^{-2}$
powers, we obtain$$
\begin{cases}
\omega_1\wedge\omega_{13}+\omega_2\wedge\omega_{23}=0&\\
\omega_1\wedge\omega_2+\omega_{13}\wedge\omega_{23}=0,&
\end{cases}$$
that is equations (3.1.3.f) and (3.1.3.g).

Equations (3.1.3.c-e) represent a particularization for $\lambda=1$ of
equations
(3.1.4.i-iii).
The coefficients of $\lambda$ and equations (3.1.5.ii,iii) determine the
expressions
of $\d\omega_1$ and $\d\omega_2$, that is the remaining conditions
($K$).\hfill$\square$

This also completes the proof of the Theorem 3.1.1.

\begin{Rem}
As frequently observed in soliton theory, the introduction of a parameter
reduces the
number of defining equations.
\end{Rem}

Theorem 3.1.1, that we just proved, is of central importance for the present
study. In section 2.3, we analyzed the pseudospherical surfaces in detail and
described a $\lambda$-transformation, $\lambda>0$, that preserves the second
fundamental form, the Gaussian curvature and the angle between asymptotic
lines. We also presented the extended normalized frame ${\cal U}^\lambda$
(2.3.15)
associated with this transformation. In this section we studied in more
detail the effects of introducing the real positive parameter $\lambda$ by the
Lie-Lorentz transformation.
We obtained a $\lambda$-family of 1-forms $\omega_i^\lambda$,
$\omega_{ij}^\lambda$,
$i<j$, which characterizes the above-mentioned $\lambda$-family
$M=(D,\psi^\lambda)$
of associated surfaces via the $\lambda$-transformation.

\subsection{The Extended Maurer-Cartan Form $\omega^\lambda$ of an Associated
Family of Pseudospherical Surfaces and the Extended Normalized
Frame ${\cal U}^\lambda$}

In section 1.1, we identified the set ${\cal F}$ of all frames with $G$,
the group of
orientation-preserving rigid motions, via a map $g^f:{\cal F}\to G$,
$g^f(x,e_1,e_2,e_3)=(x,A)$, with $x\in {\Bbb R^3}$, $A\in{\rm SO}(3)$, such
that
$e_i=A\check{e}_i$, where $F_0=\{0,\check{e}_1, \check{e}_2,\check{e}_3\}$
was a
fixed frame.

$E_i$, $E_{ij}$ are, by definition, the {\it six vector fields
dual} to the 1-forms $\omega_i$, $\omega_{ij}$, $i,j=1,2,3$, $i<j$, i.e.
the vector fields satisfying

\begin{align*}
E_j(x^f)(F)&=\d_Fx^f(E_j)=\sum_i\omega_i(E_j)e_i^f(F)=e_j^f(F)
\end{align*}

respectively
\begin{align*}
E_{ij}(e_m^f)(F)=\d_Fe_m^f(E_{ij})=\left(\sum_n\omega_{mn}(E_{ij})e_n^f(F)\right)=
&\left(\delta_i^m\delta_j^ne_n^f-\delta_i^n\delta_j^me_n^f\right)
=\left(\delta_i^me_j^f-\delta_j^me_i^f\right)=S_{ij}(e_m^f).
\end{align*}
Here $S_{ij}$ represents the
$3\times3$ matrix with $(i,j)$-entry equal 1,
$(j,i)$-entry equal to $-1$ and zero elsewhere, $i<j$. According to the way
$E_{ij}$ acts on the frame $F$, it can be identified with the matrix $S_{ij}$.

We note that the vector fields $E_i$, $E_{ij}$, $i,j=1,2,3$, $i<j$ are
invariant with respect to the particular choice of the fixed frame $F_0$.

\begin{Rem}
Reviewing, we obtained above the formulas
\begin{align}
E_j(x^f)(F)&=\d_Fx^f(E_j)=e_j^f(F),\tag{3.2.1}\\
(E_{ij}(e_m^f))_{m=1,2,3}&=\left(\d_Fe_m^f(E_{ij})\right)_{m=1,2,3}=S_{ij}
\begin{pmatrix}e_1\\e_2\\e_3\end{pmatrix}.\tag{3.2.2}
\end{align}
These equations are satisfied for every frame $F=\{x,e_1,e_2,e_3\}$.
\end{Rem}

\begin{Def} Consider the $so(3)$-valued 1-form $\omega$ given by$$
\omega=\tilde{\omega}_{12}E_{12}+\tilde\omega_{13}E_{13}+\tilde\omega_{23}E_{23}
,\eqno(3.2.3)$$
where $\tilde{\omega}_{12}$, $\tilde{\omega}_{13}$ and $\tilde{\omega}_{23}$
are given by formulas (2.3.6 e,f,g).
We will call $\omega $ the {\it Maurer-Cartan form} of the group of
Euclidean motions.
\end{Def}

As a linear combination of matrices $E_{12}, E_{13}, E_{23}$,
the form $\omega$ becomes an
$so(3)$-valued 1-form on $G$. For a vector field $Y$ on $G$, we have$$
\omega(Y)=\sum_{i<j}\tilde{\omega}_{ij}(Y)E_{ij}\eqno(3.2.4)$$

\begin{Rem}
(a) Writing $\tilde{\omega}_{12}$, $\tilde{\omega}_{13}$ and
$\tilde{\omega}_{23}$
explicitely as in (2.3.6), the Maurer-Cartan form of the group of Euclidean
motions
restricted to ${\cal F}^{\cal M}_{1}$ for a pseudospherical surface can be
written as$$
\omega=-{\cal U}^{-1}\d{\cal U}=\begin{pmatrix}
0&\varphi_x\d x&\sin\varphi\,\d y\\
-\varphi_x\d x&0&-\d x+\cos\varphi\,\d y\\
-\sin\varphi\,\d y&\d x-\cos\varphi\,\d y&0
\end{pmatrix}.\eqno(3.2.5)$$

(b) Let us recall briefly the results from the previous section.
We have proved in Lemma 3.1.1 that if $\omega_i^\lambda$,
$\omega_{ij}^\lambda$ are given by formulas (3.1.2) and conditions
(3.1.4iii--v) are
satisfied, then all the conditions ($K$) are satisfied. We have also seen that
$\omega_i^\lambda$, $\omega_{ij}^\lambda$, $\lambda>0$, given in (3.1.2)
correspond
to an associated family of surfaces that preserve the angle $\varphi$ between
asymptotic lines, the Gaussian curvature and the second fundamental form
(Theorem
2.3.1)
and that $\omega_i^\lambda$ and $\omega_{ij}^\lambda$ can be naturally induced
by a transformation $x\mapsto\lambda^{-1}x$,
$y\mapsto\lambda y$, $\lambda>0$ of the asymptotic line parametrization.
\end{Rem}

Let us now recall from Corollary 2.3.2 that the {\it extended normalized
moving frame}
${\cal U}^\lambda:D\to{\rm SO}(3)$ of
this family of one-forms $\omega_i^\lambda$, $\omega_{ij}^\lambda$, $\lambda>0$
satisfies the equations$$
\begin{cases}
({\cal U}^\lambda)^{-1}\cdot\partial_x{\cal U}^\lambda=\dis\begin{pmatrix}
0&-\varphi_x&0\\
\varphi_x&0&\lambda\\
0&-\lambda&0
\end{pmatrix}&\\
({\cal U}^\lambda)^{-1}\cdot\partial_y{\cal
U}^\lambda=\dis\frac{1}{\lambda}\begin{pmatrix}
0&0&-\sin\varphi\\
0&0&-\cos\varphi\\
\sin\varphi&\cos\varphi&0
\end{pmatrix}.&
\end{cases}\eqno(3.2.6)$$

Comparing (3.2.5) to (3.2.6), we formulate

\begin{Def}  The ${\rm so}(3)$-valued family of 1-forms$$
\omega^\lambda=-({\cal U}^\lambda)^{-1}\d{\cal U}^\lambda$$$$
=\begin{pmatrix}
0&\varphi_x\,\d x&\lambda^{-1}\sin\varphi\,\d y\\
-\varphi_x\,\d x&0&-\lambda\d x+\lambda^{-1}\cos\varphi\,\d y\\
-\lambda^{-1}\sin\varphi\,\d y&\lambda\d x-\lambda^{-1}\cos\varphi\,
\d y&0\end{pmatrix},\eqno(3.2.7)$$
is called {\it extended Maurer-Cartan form}.
\end{Def}

\begin{Prop} The system of equations {\rm(3.1.4.i-iii)} is equivalent to$$
\d\omega^\lambda+\omega^\lambda\wedge\omega^\lambda=0,\eqno(3.2.8)$$
for every $\lambda>0$.
\end{Prop}

\noindent{\it Proof.} Assume the equations {\rm(3.1.4.i--iii)} are
satisfied. Then, by Lemma 3.1.1, the system of equations
{\rm(3.1.4.i--iii)} is equivalent to the conditions (K), defined
in (3.1.3). On the other hand, {\rm(3.1.4.i--iii)} are by
definition the Gauss-Codazzi equations for a pseudospherical
surface. On the other hand, (3.2.8) can be checked directly, and
it reduces to the Gauss-Codazzi equations: e.g., the sine-Gordon
equation is recovered immediately from the (1,2) entry of the
matrix-valued form
$\d\omega^\lambda+\omega^\lambda\wedge\omega^\lambda$.
\hfill$\square$ \vskip6pt We will call formula (3.2.8) the {\it
flatness condition}, or the {\it zero-curvature condition} for the
extended Maurer-Cartan form $\omega^\lambda$.

\begin{Rem}
From equation (3.2.7), we see that the extended Maurer-Cartan form
$\omega^\lambda$ can be written in the form$$
\omega^\lambda:=\lambda^{-1}\cdot\alpha_{-1}+\alpha_0+\lambda\cdot\alpha_1,\eqno
(3.2.9)$$
where $\alpha_0\in\underline{k}=\Bbb RE_{12}$ and $\alpha_{-1},
\alpha_1\in\underline{p}=\Bbb RE_{13}+\Bbb RE_{23}$.

More precisely, we have$$
\alpha_0=\varphi_xE_{12}\d x,\eqno(3.2.10)$$
while$$
\alpha_{-1}=(\sin\varphi\cdot E_{13}+\cos\varphi\cdot E_{23})\d
y,\eqno(3.2.11)$$
and$$
\alpha_1=-E_{23}\d x.\eqno(3.2.12)$$
\end{Rem}

\section{Loop Algebras and Groups Corresponding to Pseudospherical
Surfaces}

We now examine the system (3.1.4) in the context of the loop algebra
${\rm so}(3,\Bbb R)\otimes\Bbb R[\lambda^{-1},\lambda]$. This will lead
to interpreting the extended moving frame equations in terms of loop groups,
which opens some completely new possibilities. E.g., the extended frame
$\cal U^\lambda$ can be decomposed in the form
${\cal U}={\cal U}_+\cdot V_-={\cal U}_-\cdot V_+$.
Here ${\cal U}_-$ is an element of the form
${\cal U}_-=I+\lambda^{-1}{\cal U}_{-1}+\lambda^{-2}{\cal U}_{-2}+\cdots$,
while $V_+$ is an element of the form
$V_+=V_0+\lambda V_1+\lambda^2V_2+\cdots$, respectively. Eventually, this will
allow us to find unconstrained data, ``potentials" from which all
pseudospherical
surfaces can be constructed.

\subsection{Loop Algebras and Structure Equations. Introduction}

Let $\frak a$ be a Lie algebra over $\Bbb R$ with a finite basis
$X_1,X_2,\ldots,X_m$; i.e. every $X\in\frak a$ is expressed uniquely as a
linear
combination$$
X=a_1X_1+a_2X_2+\cdots+a_mX_m,\eqno(4.1.1) $$
where $a_j\in\Bbb R$.

The structure of the Lie algebra $\frak a$ is given by {\it Lie's equations}$$
[X_i,X_j]=C_{ij}^kX_k,\eqno(4.1.2)$$
where for convenience we used the Einstein summation convention for the index
$k$, which will be used from now on.

An immediate consequence of the skew-symmetry of the Poisson bracket is the
skew-symmetry of the structural constants $C_{ij}^k$ with respect to the
indices
$i, j$.
Also, as a consequence of the Jacobi identity, the structural constants
satisfy
the following identity:
$$
C_{sj}^k C_{ir}^s + C_{si}^k C_{rj}^s + C_{sr}^k C_{ji}^s = 0 .$$

This identity appears in literature as Lie's quadratic identity.

Let $\frak a^*$ be the dual space of $\frak a$. By definition, the
dual basis of $\frak a^*$ is $\{\eta^1,\eta^2,\ldots,\eta^m\}$
such that $\eta^i(X_j)=\delta_j^i$. Also, for every $\eta\in\frak
a^*$, there is a unique linear combination$$
\eta=\beta_1\eta^1+\beta_2\eta^2+\cdots+\beta_m\eta^m.\eqno(4.1.3)$$
Let $\Lambda^p\frak a^*$ denote all the $p$-forms on $\frak a$.
Clearly, $$ \Lambda^1\frak a^* =\frak a^*.\eqno(4.1.4)$$
\begin{Def} The exterior differential $\d\eta\in\Lambda^2\frak a^*$ of a 1-form
$\eta\in\frak a^*$ is defined by the equation$$
\d\eta(X,Y)=-\eta([X,Y]),\eqno(4.1.5)$$
where $X,Y\in\frak a$.
\end{Def}

Equation (4.1.5) is equivalent to {\it Cartan's structure equations}:$$
\d\eta^k+\frac{1}{2}C_{ij}^k\eta^i\wedge\eta^j=0.\eqno(4.1.6a)$$
This equivalence is straightforward and is presented in classical
texts (e.g., [Ca], p.45).
In (4.1.6a), $\eta^i\wedge\eta^j$ represents the exterior product of the
1-forms $\eta^i$ and $\eta^j$.

It is easy to see that (4.1.6a) can be rewritten as$$
\d\eta^k+C_{ij}^k\eta^i\wedge\eta^j=0, \eqno(4.1.6b)$$
where $i<j$.

Multiplying equation (4.1.6b) by $X_k$ and taking into account Lie's
equations (4.1.2), we obtain
$$
X_k\cdot d\eta^k+[X_i,X_j]\eta^i\wedge\eta^j=0, \qquad i<j,$$
which can be rewritten as
$$
\d\eta+\frac{1}{2}[\eta\wedge\eta]=0,  \eqno(4.1.7)$$
where
$$
\eta=X_1\eta^1+X_2\eta^2+\cdots+X_m\eta^m.  \eqno(4.1.8) $$

\begin{Rem} If the basis $\{\eta^1,\eta^2,\ldots,\eta^m\}$
 of $\frak a^*$ is divided into two groups distinguished by indices
 $i,j,k\in N_1$ and $\alpha,\beta,\gamma\in N_2$ respectively,
 then the structure equations become$$
\begin{cases}
\d\eta^k+C_{ij}^k\eta^i\wedge\eta^j+C_{i\beta}^k\eta^i\wedge\eta^\beta
+C_{\alpha\beta}^k\eta^\alpha\wedge\eta^\beta=0,  i<j, \alpha<\beta &\\
\d\eta^\gamma+C_{ij}^\gamma\eta^i\wedge\eta^j+C_{i\beta}^\gamma\eta^i\wedge\eta^
\beta
+C_{\alpha\beta}^\gamma\eta^\alpha\wedge\eta^\beta=0,  i<j, \alpha<\beta &
\end{cases}.\eqno(4.1.9)$$

Note that the restriction $\eta^\gamma=0$, for every $\gamma\in N_2$,
defines a linear subspace of $\frak a$.
\end{Rem}

\begin{Exam}
Consider the group of Euclidean motions $T$ given by the structure equations
(1.1.17) and introduce the restrictions
$\omega_{12}=\omega_{13}=\omega_{23}=0$,
which define the normal subgroup of all translations.
The groups of indices specified in Remark 4.1.1 are
$1,2,3\in N_1$ and $12, 13, 23\in N_2$ respectively, where we replaced
 $\eta$ by $\omega$.

Let us consider the quotient group $G/T={\rm O}(3,\Bbb R)$
of the Euclidean motion group modulo the group of translations.
Thus, in the second group of equations of the system (4.1.9), the terms
containing $\omega_j, j=1,2,3$ disappear, and the equations become$$
\d\omega_{ij}=\omega_{ik}\wedge\omega_{kj},$$
with Einstein summation with respect to $k$ and $i,j = 1,2,3, i<j$.

This gives a concrete illustration of the structure equations (1.1.24 c,d,e).
\end{Exam}

The form (4.1.8) for the Euclidean motion group is written here as$$
\hat\omega=\omega_1E_1+\omega_2E_2+\omega_3E_3+\omega_{12}E_{12}+\omega_{13}E_{1
3}+\omega_{23}E_{23},\eqno(4.1.10)$$
\
The form $\hat\omega$ is sometimes called the {\it total Maurer-Cartan form}.

\subsection{The Loop Algebra Setting}

Let now $\frak b$ represent a Lie algebra with basis $X_1,X_2,\ldots,X_m$
satisfying $[X_i,X_j]$ $=C_{ij}^kX_k$. This is equivalent to the structure
equations (4.1.6).

\begin{Def} The {\it polynomial loop algebra} $\frak a=\frak b\otimes
{\Bbb R}[\lambda^{-1},\lambda]$
is the Lie algebra with basis $X_{k,t}=X_k\lambda^t$, $k=1,2,\ldots,m$,
$t=0,\pm1,\pm2,\ldots$, where $\lambda$ is a formal parameter.
\end{Def}

This basis satisfies the Lie equations$$
[X_{i,r},X_{j,s}]=C_{ij}^kX_{k,r+s}.\eqno(4.2.1)$$
The notation ${\Bbb R}[\lambda^{-1},\lambda]$ used above represents the
ring of Laurent
polynomials in the variable $\lambda$ over the field ${\Bbb R}$. Let
$\{\eta^{i,r}\}$
represent the basis of 1-forms dual to the basis $\{X_{i,r}\}$. Then,
analogous to the
derivation of (4.1.6b) we obtain, as a consequence of (4.2.1),
the structure equations of the loop algebra $\frak a$$$
\d\eta^{k,t}+\sum_{r+s=t,
i<j}C_{ij}^k\eta^{i,r}\wedge\eta^{j,s}=0.\eqno(4.2.2)$$

Multiplying these equations by $\lambda^t=\lambda^{r+s}$, we obtain$$
\d\eta^{k,t}\lambda^t+\sum_{r+s=t,
i<j}C_{ij}^k\eta^{i,r}\lambda^r\wedge\eta^{j,s}\lambda^s
=0.
\eqno(4.2.3)$$

That is, the structure equations of the form (4.1.6), where$$
\eta^k=\sum_{t=-\infty}^\infty\eta^{k,t}\lambda^t$$
represent infinite Laurent series in the variable $\lambda$ with 1-forms
$\eta^{k,t}$ as coefficients.

Let us now consider the particular case of $\frak b={\rm so}(3,\Bbb R)$,
so that
$\frak a={\rm so}(3,\Bbb R)\otimes\Bbb R[\lambda^{-1},\lambda]$. The main
reason why
we focus on this loop algebra is provided by the extended Maurer-Cartan form
$\omega^\lambda$ of a pseudospherical surface, introduced in (3.2.7).
Moreover, we
shall introduce the twisted loop algebra$$
\Lambda{\rm so}(3)_P^{\rm alg}=\{X\in{\rm so}(3)\otimes\Bbb
R[\lambda,\lambda^{-1}];\ X(-\lambda)=PX(\lambda)P^{-1}\},\eqno(4.2.4)$$
where$$
P={\rm diag}\{1,1,-1\}.$$

Note that $P^{-1}=P$ and $$
PE_{12}P=E_{12}, PE_{13}P=-E_{13}, PE_{23}P=-E_{23}.\eqno(4.2.5)$$

 From (3.2.7), it is easy to see that
$\omega^\lambda(-\lambda)=P\cdot\omega^\lambda(\lambda)\cdot P^{-1}$ holds.
Hence, $\omega^\lambda\in\Lambda{\rm so}(3)_P^{\rm alg}$.

It will be convenient to use certain Banach completions of the Lie algebra
(4.2.4). For this purpose, for a matrix $A\in {\rm so}(3,\Bbb R)$
independent of $\lambda$,
we introduce the norm
$$
\quad\|A\|=\max_{i} \{ {\sum_{j=1}^3}|A_{ij}| \},\eqno(4.2.6)$$

where $A_{ij}$ denotes the $(i,j)$-coefficient of $A$.

It can be checked by a direct computation that $$
\quad\|A B\|\leq \|A\| \cdot \|B\|. \qquad\|I\|=1.$$

Further, if ${X(\lambda)}=\sum_{k\in\Bbb Z}X_k\cdot \lambda^k$, we define
its norm as
follows:$$
\|{X(\lambda)}\|=\sum_{k\in\Bbb Z}\|X_k\|<\infty .\eqno(4.2.7)$$

\begin{Rem}

The norm defined by (4.2.7) can be also introduced as follows:

We start by defining the norm of a real-valued function in $\lambda $,$$
\|h\|:=\sum_{k\in\Bbb Z}|h_k|<\infty ,\qquad h(\lambda )= {\sum_{k\in\Bbb
Z}}{h_k}{\lambda^k }.$$

Then we define the norm of the matrix-valued function $X(\lambda )$ as
$$
\quad\|X\|=\max_{i} \{ {\sum_{j=1}^3}\|{X_{ij}(\lambda)}\|  \}.$$

It is easy to see that we obtain this way the same norm as in (4.2.7).

\end{Rem}

Note that in (4.2.6) and (4.2.7), by abuse of notation, we use the same
symbol $\|\cdot\|$
for the following three different items: norm of a function, norm of a
$\lambda $-independent
matrix and norm of $X(\lambda )$. It will always be clear from the context
which norm we mean.

We set$$
\Lambda{\rm so}(3)_P:=\text{completion of }\Lambda{\rm so}(3)_P^{\rm alg}\text{
relative to }\|\cdot\|.\eqno(4.2.8)$$

\begin{Prop}
$\Lambda{\rm so}(3)_P$ is a Banach Lie algebra.
\end{Prop}

\noindent{\it Proof.} We can define the norm (4.2.7) for arbitrary
matrices
in ${\rm gl}(3)\otimes\Bbb R[\lambda,\lambda^{-1}]$.

The fixed point algebra of the automorphism $X(\lambda)\mapsto P\cdot
X(-\lambda)\cdot P^{-1}$ of $\Lambda{\rm GL}(3,\Bbb R)$ is an associative
Banach
subalgebra. Inside the connected component of the Banach Lie group of
invertible
elements of this fixed point algebra, we consider the connected component
of the group$$
\Lambda{\rm SO}(3)_P=\{g\in\Lambda{\rm SO}(3,\Bbb R);\
Pg(\lambda)P^{-1}=g(-\lambda)\}.\eqno(4.2.9)$$
From [Ha,Ka], it follows that $\Lambda{\rm SO}(3)_P$ is a Banach Lie group
with Lie algebra$$
{\rm Lie}\,\Lambda{\rm SO}(3)_P=\Lambda{\rm so}(3)_P.\eqno(4.2.10)$$

\begin{Rem} If $M=(D,\psi)$ is, as usual, a pseudospherical surface
given by the Chebyshev immersion $\psi:D\to\Bbb R^3$, where $D$ is
a simply connected domain, then there exists a normal $N:D\to S^2$
along $\psi$ and a frame ${\cal U}:D\to{\rm SO}(3)$ along $\psi$
such that $e_3=N$ denotes the Gauss map of $\psi$:

\noindent $\pi$ above denotes the canonical projection relative to the base
point $e_3$.
Thus, $S^2\cong{\rm SO}(3)/K$.
Note that the Lie algebra of the group $K \simeq SO(2)$ is
Lie $K=\underline{k}=\Bbb RE_{12}$.
\end{Rem}

\begin{Rem} As we pointed out, giving an extended Maurer-Cartan form
$\omega^\lambda$
satisfying the flatness condition is equivalent to giving the forms
$\omega^\lambda_i$, $\omega_{ij}^\lambda$, $i<j$ satisfying the conditions
($K$),
which is also equivalent to giving a family of surfaces $M_{\lambda }$ of
constant negative Gaussian curvature $K=-1$. To such an associated
family of surfaces, we
attached ( see (3.2.7) ) the extended frame ${\cal U}^\lambda:D\times\Bbb
R_+\to\Lambda{\rm SO}(3)_P$ satisfying $({\cal U}^\lambda)^{-1}\d{\cal
U}^\lambda+\omega^\lambda=0$, where $\Bbb R_+$ represents the set of strictly
positive real numbers $\lambda$. It will be convenient for our purposes to
fix a base point $x_0\in D$ , e.g. $x_0=(0,0)$, and require that the frame
satisfies the ``initial condition" $$
{\cal U}(x_0,\lambda)=I,\eqno(4.2.11)$$
for every $\lambda$. We will use this assumption from now on.
\end{Rem}

\begin{Rem} The subalgebra $\Lambda{\rm so}(3)^{\rm alg}_P$ of ${\rm
so}(3)\otimes\Bbb R[\lambda,\lambda^{-1}]$ defined by (4.2.4) can also be
characterized as the subalgebra consisting of elements with the following

\begin{itemize}
\item {\sf Property} $\cal P$: In a representation relative to the basis
      $E_{12},E_{13},E_{23}$, the coefficient of $E_{12}$ is an even
      function of $\lambda$, while the coefficients of $E_{13}$ and
      $E_{23}$ are odd functions of $\lambda$.
\end{itemize}
\end{Rem}

\subsection{Loop Groups and Group Splittings Used for Pseudospherical
Surfaces}

In order to carry out the DPW method in the context of pseudospherical
surfaces, we
introduce the following subalgebras of $\Lambda{\rm so}(3)_P$:
\begin{align}
\Lambda^+{\rm so}(3)_P&=\{X(\lambda)\in\Lambda{\rm so}(3)_P;\ X(\lambda)\text{
contains only non-negative}\nonumber\\
&\hskip9.5pc\text{powers of }\lambda\}\tag{4.3.1}\\
\Lambda^-{\rm so}(3)_P&=\{X(\lambda)\in\Lambda{\rm so}(3)_P;\ X(\lambda)\text{
contains only non-positive}\nonumber\\
&\hskip9.5pc\text{powers of }\lambda\}\tag{4.3.2}\\
\Lambda^-_*{\rm so}(3)_P&=\{X(\lambda)\in\Lambda^-{\rm so}(3)_P;\
X(\infty )=0\}\tag{4.3.3}
\end{align}

The connected Banach loop groups whose Lie algebras are described by
definitions
(4.3.1--4.3.3) are denoted, respectively, $\Lambda^+{SO}(3)_P$,
$\Lambda^-{SO}(3)_P$ and $\Lambda^-_*{SO}(3)_P$.

A first question arises when we aim to split \`a la Birkhoff elements from
$\Lambda{SO}(3)_P$ with $\lambda\in\Bbb R_+$ instead of $\lambda\in
S^1$. The classical factorization theorem is stated and proved in [Pr, Se] for
smooth loops on $S^1$ and reformulated in [DPW], [DGS] for a complexified
Banach loop group $G^C$.

For our applications, the relevant part is

\begin{Th} {\rm[ DPW; Thm. 2.2.], [Pr, Se; Thm. 8.1.1--8.1.2]:}
\vskip6pt

Let $G$ be a compact Lie group. Then the multiplication
$\Lambda_*^-G^C\times\Lambda^+G^C\to\Lambda G^C$ is an analytic diffeomorphism
onto the open and dense subset $\Lambda_*^-G^C\cdot\Lambda^+G^C$,
called the ``big cell". In particular, if $g\in\Lambda G^C$ is contained
in the big cell, then $g$ has a unique decomposition$$
g=g_-g_+,\eqno(4.3.4)$$
where $g_-\in\Lambda_*^-G^C$ and $g_+\in\Lambda^+G^C$. The analogous
result holds for the multiplication map
$\Lambda_*^+G^C\times\Lambda^-G^C\to\Lambda G^C$.
\vskip6pt
The results stated above hold in particular for $G={SO}(3)$. The splitting
(4.3.4)
is called the Birkhoff factorization of $\Lambda G^C$.
\end{Th}

\noindent{\bf Remark A.} Regarding the $\lambda\in S^1$ versus
$\lambda\in\Bbb R_+$ issue, our Appendix contains the proof of the
fact that the splitting works also for some specific ``loop" group
with real, positive $\lambda$.

Let $\tilde\Lambda SO(3)_P$ be the subset of $\Lambda SO(3)_P$
whose elements, as maps defined on $\Bbb R_+$, admit an analytic extension to
$\Bbb C_*$. It is easy to see that $\tilde\Lambda SO(3)_P$ is a
subgroup of $\Lambda SO(3)_P$.
We have the following result:

\begin{Th}
\vskip6pt
$\tilde\Lambda_*^-SO(3)_P\times\tilde\Lambda^+SO(3)_P\to\tilde\Lambda SO(3)_P$
is a diffeomorphism onto the open and dense subset
$\tilde\Lambda_*^-SO(3)_P\cdot\tilde\Lambda^+SO(3)_P$, called the ``big
cell". In particular,
if $g\in\tilde\Lambda SO(3)_P$ is contained
in the big cell, then $g$ has a unique decomposition$$
g=g_-g_+,\eqno(4.3.5)$$
where $g_-\in\tilde\Lambda_*^-SO(3)_P$ and $g_+\in\tilde\Lambda^+SO(3)_P$.
The analogous
result holds for the multiplication map
$\tilde\Lambda_*^+SO(3)_P\times\tilde\Lambda^-SO(3)_P\to\tilde\Lambda SO(3)_P$.
\vskip6pt
\end{Th}

{\it Proof. See Appendix}.
\hfill$\square$

Remark that any extended frame $\cal U^\lambda$, as a function of
the real positive parameter $\lambda$, admits an analytic extension
to $\Bbb C_*$. This is straight-forward and is stated and proved in
Lemma A.1.

Hence, any extended frame ${\cal U}(x,y,\lambda)$ from the ``big cell"
of $\tilde\Lambda SO(3)_P$ can be split as
$$
{\cal U}={\cal U}_+\cdot V_-={\cal U}_-\cdot V_+.\eqno(4.3.6)$$

Here
${\cal U}_-$ is an element of the form
${\cal U}_-=I+\lambda^{-1}{\cal U}_{-1}+\lambda^{-2}{\cal U}_{-2}+\cdots$,
while $V_+$ is an element of the form
$V_+=V_0+\lambda V_1+\lambda^2V_2+\cdots$, respectively. Analogous expressions
can be written for ${\cal U}_+$ and $V_-$, respectively. Namely,
${\cal U}_+$ is an element of the form
${\cal U}_+=I+\lambda {\cal U}_{1}+\lambda^{2} {\cal U}_{2}+\cdots$,
while $V_-$ is an element of the form
$V_-=V_0+\lambda^{-1}V_{-1}+\lambda^{-2}V_{-2}+\cdots$.

\section{Harmonic Maps and Generalized Weierstrass Data}

In this section we present the notion of harmonic map from a
pseudospherical surface
$M$ to $S^2$. This is a particular case of a harmonic map from a
pseudo-Riemannian
manifold to another pseudo-Riemannian manifold, i.e. a differentiable map
whose
tension field vanishes (see [EL]). The Gauss maps of certain classes of
surfaces
(e.g. constant mean
curvature, minimal, constant Gaussian curvature) are harmonic with respect
to some
suitable (pseudo)metrics. It was proved that the harmonic maps from these
classes of
surfaces to $S^2$ are in one-to-one correspondence
with the equivalence classes of flat extended forms $\omega^\lambda$ (3.2.8)
under the action of a gauge group. In connection with Sections 3 and 4,
this is a
strong motivation for studying such harmonic maps.

\subsection{Harmonic Maps}

\begin{Def}
Let $(M,g)$ and $(\tilde M, \tilde g)$ be pseudo-Riemannian manifolds.
A harmonic map $f:M\to \tilde M$ is a differentiable map such that
its {\it tension field} $\tau(f)$ vanishes:$$
\tau(f):={\rm Trace}(\nabla\d f)=0,\eqno(5.1.1)$$
where $\nabla$ is the Levi-Civita connection on the vector bundle
$T^*(M)\otimes f^*(T \tilde M)$, provided with the natural pseudo-metric
induced by $g$
and $\tilde g$.
\end{Def}

For Riemannian manifolds, the system (5.1.1) is elliptic. This property is
not maintained
on pseudo-Riemannian manifolds. In this case, harmonic maps are sometimes
called
{\it pseudo-harmonic}.

The notion of harmonic map was first introduced by Eells and Sampson for
Riemannian
manifolds, then generalized to pseudo-Riemannian manifolds by Eells and
Lemaire ([EL])
and then studied by several authors (e.g., [GU], [Me, St, 1]).

If $(M,g)$ and $(\tilde M,\tilde g)$ are two Riemannian manifolds, $df(x)$
represents the differential of $f$ ( linear map from $TM$ to $T \tilde M$ at a
point $x$ of $M$), while its {\it tension field} is
$$
\tau(f)= {\rm div}( d f)= g^{ij} (\nabla(d f))_{ij}.\eqno(5.1.2)$$

Here we used again the Einstein summation convention with respect to both
indices $i,j$.
$g^{ij}$ are the entries of the inverse $g^{-1}$ of the matrix $g$.


The integral over $M$ of the energy density $|d f|^2$ with respect to the
area element
on $M$ is frequently called {\it energy functional}. Equation (5.1.1)
arises as the
Euler-Lagrange equation for the variational problem of the energy integral.
Harmonic maps $f$ represent critical points of the energy functional.

We shall now introduce a concept which is actually equivalent to the one of
 extended Maurer-Cartan form $\omega^{\lambda}$.

\begin{Rem}
The following represents a necessary and sufficient condition for a map to
be harmonic
([UR]):

{\sf Lemma :}

Let $f$ be a smooth map from a pseudo-Riemannian manifold to the sphere $S^n$.
Then $f$ is harmonic {\it iff}$$
\Delta f = \rho\cdot f,\eqno(5.1.3)$$
for some function $\rho$, where $\Delta$ represents the Lorentz-Laplace
operator.

In this case, $\rho=e(f)=|d f|^2$ is the energy density of $f$.

For the case $f:M \to S^2$, where $M$ is a 2-dimensional manifold, see also
[Me, St, 1],
Prop. 1.1. Moreover, harmonicity is invariant under conformal transformations.

\end{Rem}

\begin{Rem}
A classically known fact is the following:

If $M$ is a weakly regular surface with $K<0$, then $M$, endowed with its
second fundamental form ${\rm II}$ (2.1.2) in asymptotic coordinates, is a
{\it Lorentzian}
2-manifold $(M,{\rm II})$.

Moreover, the Gauss map $N:(M,{\rm II})\to S^2$ is harmonic {\it
iff} $K=$ constant. With respect to the second fundamental form,
(5.1.3) is written as $$ N_{xy} = \rho\cdot N.\eqno(5.1.4)$$ In
this sense, the Gauss map of every pseudospherical surface is
harmonic.

This property of pseudospherical surfaces is sometimes called
Lorentz-harmonicity.

\end{Rem}

\begin{Def}
Let us consider an $so(3)$-valued form $\omega$.

Recall from the previous section the Lie algebras $\underline{k}=\Bbb
RE_{12}$ and
$\underline{p}=\Bbb RE_{13}+\Bbb RE_{23}$.

Let $\eta=\eta_0 + \eta_1$
be the Cartan decomposition of $\eta$ into its $\underline{k}$-part $\eta_0$,
respectively its $\underline{p}$- part, $\eta_1$. Then $\eta$ is called an
{\it admissible connection}  if it satisfies the following pair of equations
(sometimes called {\it Yang-Mills-Higgs} equations):$$
\d\eta+\eta\wedge\eta=\d\eta+\frac{1}{2}[\eta\wedge\eta]=0,\eqno(5.1.5)$$$$
\d(\ast\eta_1)+[\eta_0\wedge\ast\eta_1]=0.\eqno(5.1.6)$$

\end{Def}

For (5.1.5) and (5.1.6), see [Gu, Oh].

From the Remark 4.2.1, the smooth Gauss map $N$ has the frame
${\cal U}$ as a lift. It follows (e.g, [Bo 2]) that the maps $N$
and ${\cal U}$ are related by the identification $$ N\equiv{\cal
U}\cdot E_{12}\cdot {\cal U}^{-1}.\eqno(5.1.7) $$

{\sf Note}:  In (5.1.7), [Bo2] uses $-i \sigma_3$ instead of our
$E_{12}$. $\sigma_3$ is the third Pauli matrix (6.4.1). This fact
is explained by the (spinor representation) isomorphim between
$su(2)$ and $so(3)$, which is presented in Section 6.4.

A very important result obtained by A. Sym ([Sy]) allows us to obtain the
immersion once
we have the expression of the extended frame. This is presented in several
papers,
including for the particular case of pseudospherical surfaces (e.g. [1, Me,
St],
[Bo, Pi]) and can be stated as follows:

\begin{Th}
Starting from a given $\varphi (x,y)$, a solution to the sine-Gordon
equation, let us
consider the initial value problem consisting of the Lax system (2.3.15)
together
with the initial condition ${\cal U}(0,0,\lambda) = I$. Let $\cal
U(\lambda)$ be the
solution to this initial value problem. Then $\cal U(\lambda)$ represents
the extended frame
corresponding to the Chebyshev immersion
$$
\psi^{\lambda} = {\d\over\d t}{\cal U}^{\lambda}\cdot ({\cal
U}^{\lambda})^{-1},\eqno(5.1.8)$$
where $\lambda = e^t$.
\end{Th}

By Theorem 5.1.1, once we have the extended frame, we can reconstruct the
surface.
Also, the relationship between the extended frame $\cal U$ and the Gauss
map $N$ is clear, via (5.1.7). So in a sense we could reconstruct
everything starting from
the Gauss map. However, there is a freedom in the frame given by a gauge
action.

\begin{Def} Let us consider a rotation of angle $\theta$ around $e_3$,$$
R=\begin{pmatrix}
\cos\theta&\sin\theta&0\\
-\sin\theta&\cos\theta&0\\
0&0&1\end{pmatrix}.$$
\end{Def}

The rotation $R$, thought of as an element of $SO(2)$, acts on the frame
$\cal U$,
and produces the so called {\it gauged frame} $\hat{\cal U}$ of the
pseudospherical
surface $M$, via the rule$$
\hat{\cal U}={\cal U}\cdot R^{-1}.\eqno(5.1.9)$$
As a consequence of this action by a rotation matrix on the frame, the
Maurer-Cartan form
$\omega$ changes accordingly, to a $\hat{\omega}$. On the other hand, the
Gauss map
$N = {\cal U}\cdot E_{12}\cdot {\cal U}^{-1}$
from equation (5.1.7) is obviously invariant under such a gauge transformation.

The following very important result is a particular case of [Me, St, 1],
Prop.1.4.

\begin{Prop}
There is a one-to-one correspondence between the space of harmonic maps
from the Lorentzian surface $M$ to $S^2$ and the equivalence classes
of admissible connections, under the action of the gauge group introduced
by (5.1.9).
\end{Prop}

\begin{Rem}
On the other hand, every admissible connection ``$\omega$ corresponds
to its associated loop $\omega^{\lambda}$ satisfying the flatness condition
(3.2.8):$$
\d\omega^\lambda+\omega^\lambda\wedge\omega^\lambda=0.$$
Recall that we called $\omega^{\lambda}$ extended Maurer-Cartan form.
\end{Rem}

The result above provides a strong interest in harmonic maps.
Summarizing, the Gauss map of a pseudospherical surface has the following
properties:

\begin{Th} {\rm[Bo2, Prop. 7]} The Gauss map $N:M\to S^2$ of a
surface with $K=-1$ is {\rm Lorentz-harmonic}, i.e.,$$
N_{xy}=qN,\qquad q:M\to\Bbb R.\eqno(5.1.10)$$

Moreover, $N$ forms in $S^2$ the same kind of Chebyshev net as the immersion
function
does in $\Bbb R^3$:$$
|N_x|=A,\quad|N_y|=B,\qquad\text{where }A=|\psi_x|,\ B=|\psi_y|.\eqno
(5.1.11)$$
\end{Th}

\ {\it Proof.} A lengthy but straight-forward calculation using
formulas (5.1.7) and (5.1.8) leads to formulas (5.1.10, 5.1.11).
\hfill$\square$

Via Proposition 5.1.1 and Theorem 5.1.2, we state the following:

\begin{Rem}

As a consequence of the previous results and remarks, we conclude:

A smooth map $N: D\to S^2$ is Lorentz-harmonic if and only if there is an
extended
frame ${\cal U}:D\to\Lambda{\rm SO}(3)_P$ such that $\pi\circ{\cal
U}^\lambda|_{\lambda=1}=N$,
and such that $$
\omega^\lambda:=-({\cal U}^\lambda)^{-1}\d{\cal U}^\lambda\eqno (5.1.12)$$
satisfies the flatness condition (3.2.8).

Here we denoted by $\pi:{\rm SO}(3)\to{\rm SO}(3)/K$ the canonical projection,
and $K$ a Lie subgroup isomorphic to {\rm SO}(2), which is the isotropy
group of the
action of {\rm SO}(3) on the vector $e_3$ in $\Bbb R^3$.

\end{Rem}

Let $O$ be the point corresponding to $x=0$, $y=0$ in $M$. We consider the
extended frame
corresponding to the frame ${\cal U}$ the solution ${\cal U}^\lambda$ of
equation (5.2.6) that satisfies the additional initial condition$$
{\cal U}^\lambda(0,0,\lambda)={\cal U}(0,0)=I,\eqno (5.1.13)$$

where ${\cal U}$ is the frame of $N:D\to S^2$, $N(0,0)=e K$, such that
Lie $K=\underline{k}=\Bbb RE_{12}$.
Clearly, ${\cal U}^\lambda(x,y,1)={\cal U}(x,y)$.

Let us now consider the Cartan decomposition
$\underline{g}=\underline{k}+\underline{p}$
where $\underline{k}=\Bbb RE_{12}$ and $\underline{p}=\Bbb RE_{13}+\Bbb
RE_{23}$.
Let $\omega^\lambda$ be a 1-form that satisfies the flatness condition (3.2.8).

Via the Cartan decomposition above, $\omega^\lambda$ can be written in the
form$$
\omega^\lambda:=\alpha_0+\omega_1^{\lambda },\eqno(5.1.14)$$
where $\alpha_0\in\underline{k}$ and
$\omega_1^{\lambda}=\lambda^{-1}\cdot\alpha_{-1}+\lambda\cdot\alpha_1\in\underline{p}$.

As a consequence of Theorem 5.1.1, we obtain:

\begin{Prop}
Let ${\cal U}^\lambda:D\to\Lambda{\rm SO}(3)_P$ be any map such that
$({\cal U}^\lambda)^{-1}\d{\cal U}^\lambda$ is of the form (5.1.14) and
satisfies the
flatness
condition (3.2.8). Then ${\cal U}^\lambda$ represents an extended
normalized frame
corresponding to the associated family of Chebyshev immersions$$
\psi^{\lambda} = {\d\over\d t}{\cal U}^{\lambda}({\cal U}^{\lambda})^{-1},\eqno(5.1.15)$$

\end{Prop}

\subsection{The Weierstrass-type Representation}

$\mathbf A$.{\it Generalized Weierstrass Representation of Constant Mean
Curvature Surfaces}.

In [DPW], the authors have introduced a Weierstrass type representation through
which every harmonic map from a Riemann surface $M$ to an arbitrary compact
symmetric space $G/K$ is described by a ${\rm Lie}\,{\rm G}^C$ - valued
meromorphic
differential on the universal covering of $M$. In [Do, Ha], the authors
present the case of a
constant mean curvature surface $M$ in $\Bbb R^3$, parametrized in
conformal coordinates,
obtaining the above-mentioned  differential explicitely.

For the case of $G={\rm SO}(3)$ and $K={\rm SO}(2)$, $G/K\cong S^2$, this
procedure is based on
introducing the extended normalized frame ${\cal U}^\lambda:D\to\Lambda{\rm
SO}(3)_P$,
which for $\lambda=1$ represents the normalized moving frame.
In this case, the ${\rm so}(3,\Bbb C)$-valued meromorphic differential is
characterized by two
different meromorphic functions.
The poles of the above mentioned meromorphic functions are situated at
points where
the Birkhoff loop group factorization ${\cal U}={\cal U}_-V_+$ fails to exist.

The
Weierstrass-type data is expressed via a Lie algebra-valued differential form$$
\xi= {\cal U}_-^{-1}\d{\cal U}_-=\lambda^{-1}\eta.\eqno (5.2.1)$$

\begin{Def} The forms $\eta$ and $\xi$ given by equation (5.2.1)
are called ( see also [Wu2] and [DoHa]) {\it normalized}, and respectively
{\it meromorphic} potentials.
\end{Def}

Starting from the normalized potential, we can construct the
associated family of CMC surfaces $M_\lambda=(D,\psi_\lambda)$.

An analogous result is presented in [DPT] for minimal surfaces in $\Bbb R^3$,
parametrized in conformal coordinates.

$\mathbf B$.{\it Generalized Weierstrass Representation of Pseudospherical
Surfaces}

 The aim of Sections 5 and 6 is to present the analogue of the
DPW method explained above for the case of pseudospherical
surfaces. The main result of the Section 6 is the Weierstrass-type
data for pseudospherical surfaces. In Section 6 we define the
generalized Weierstrass representation as a pair of Lie
algebra-valued differential forms
$$
\xi^x= - {\cal U}_+^{-1}\d{\cal U}_+=\lambda\eta^x,\eqno (5.2.2a)$$$$
\xi^y= - {\cal U}_-^{-1}\d{\cal U}_-=\lambda^{-1}\eta^y.\eqno (5.2.2b)$$

\begin{Def} The forms $\eta^x$ and $\eta^y$ given by equations (5.2.2a, 5.2.2b)
are called {\it normalized} $x$-potential and $y$-potential, respectively.
\end{Def}

Starting from such a  pair of normalized potentials, we can
construct the associated family of pseudospherical surfaces
$M_\lambda=(D,\psi_\lambda)$.

\section{Explicit Forms of the Normalized Potentials of Pseudospherical
Surfaces}

\subsection{Normalized Potential for CMC Surfaces Revisited}

For constant mean curvature surfaces $M=(D,\psi)$ parametrized in conformal
coordinates with metric $\d s^2=4e^{2\omega(z,\bar{z})}\d z\d\bar{z}$,
Theorem 2.1,
[Wu2], offers a simple method to calculate the normalized potential.

Namely, if the Maurer-Cartan form is
\begin{align}
{\cal U}^{-1}\d{\cal
U}&=\alpha_{-1}\lambda^{-1}+\alpha_0+\alpha_1\cdot\lambda\nonumber\\
\alpha_0&=\alpha'_0\d z+\alpha_0''\d\bar{z},\tag{6.1.1}
\end{align}
we denote by $\beta_0(z)$ and $\beta_1(z)$, respectively, the holomorphic
part $\alpha_0'(z,0)\d z$ of $\alpha_0'\,\d z$ and the holomorphic part
$\alpha_{-1}(z,0)$  of $\alpha_{-1}$. Recall that the holomorphic part
of a function $f(z,\bar{z})=\sum_{k,l} a_{kl} z^{k} {\bar z}^l$ is
$f(z,0)$.

Then the following theorem will provide the normalized potential
$\eta$:
\noindent{\bf Theorem 6.1.1} {\rm(2.1, [Wu2])} {\it The
normalized potential $\eta$ of the surface, with the origin $z=0$
as the reference point, is given by$$
\eta(z)=\psi_0(z)\cdot\beta_1(z)\cdot\psi_0(z)^{-1},\eqno(6.1.2)$$
where $\psi_0$ is the solution to$$
\psi_0(z)^{-1}\d\psi_0(z)=\beta_0(z),\quad\psi_0(0)={\cal
U}(0),\eqno(6.1.3)$$ and ${\cal U}$ is the normalized frame at the
origin.}

For CMC surfaces (see, for example [Wu2], formula (3.18)) the normalized
potential
is of the form$$
P(z)=\begin{pmatrix}
0&0&-b(z)\\
0&0&-c(z)\\
b(z)&c(z)&0
\end{pmatrix}\d z,\eqno(6.1.4)$$
where
\begin{align}
b(z)&=\frac{1}{2}\left(e^{2\xi(z)-\xi(0)}+Q(z)e^{\xi(0)-2\xi(z)}\right),\nonumber\\
c(z)&=\frac{i}{2}\left(-e^{2\xi(z)-\xi(0)}+Q(z)e^{\xi(0)-2\xi(z)}\right),\tag{6.1.5}
\end{align}

and $\xi$ represents the {\it holomorphic part} $\omega(z,0)$ of
$\omega(z,\bar{z})$, where $$
\d s^2=4e^{2\omega(z,\bar{z})}\d z\d\bar{z}$$
represents the metric of the surface, while $Q(z)=(N,\psi_{zz})$ is the
(holomorphic)
coefficient of the Hopf differential $Q(z)(\d z)^2$.

Equivalently, under the adjoint map ${\rm Ad}:{\rm SU}(2)\to{\rm SO}(3)$
(see [Wu2],
Remark 3.22, and [DoHa]), via a lifting to ${\rm SU}(2)$, the normalized
potential
can be written as$$
\eta(z)=\frac{1}{2}\begin{pmatrix}
0&e^{2\xi(z)-\xi(0)}\\
-Q(z)e^{\xi(0)-2\xi(z)}&0
\end{pmatrix}\d z.\eqno(6.1.6)$$
In the following subsection we shall state and prove a similar result for
pseudospherical
surfaces parametrized in asymptotic line coordinates.

\subsection {Normalized $x$- and $y$- Potentials for Pseudospherical
Surfaces. Ordinary Differential Systems Associated with Normalized
Potentials}

 By analogy with the normalized potential introduced for constant mean
curvature surfaces,
it becomes natural to consider a normalized potential for other classes of
surfaces
whose Gauss map is harmonic, as a map between pseudo-Riemannian surfaces,
in particular
for the class of pseudospherical surfaces.

We will introduce the generalized Weierstrass representation for
pseudospherical surfaces in
a Chebyshev parametrization, as two normalized potentials:

$\eta^x$ and $\eta^y$, where $\eta^x$ does not depend on $y$, and $\eta^y$
does not depend
on $x$.

Theorem 6.2.1 below will make this explicit.

Theorems 6.3.1 and 6.3.2 in the next section will give explicit formulas
for the normalized
potentials. They are consequences of Theorem 6.2.1.

In our case, the group $K$ represents the group of rotations
around $e_3$, isomorphic to $SO(2)$,$$ K=\left\{\begin{pmatrix}
\cos r&-\sin r&0\\
\sin r&\cos r&0\\
0&0&1\end{pmatrix};\ r\in[0;2\pi)\right\}.\eqno(6.2.1)$$

Its Lie algebra ${\rm Lie} K$ is$$
\underline{k}=\left\{\begin{pmatrix}
0&a&0\\
-a&0&0\\
0&0&0\end{pmatrix};\ a\in\Bbb R\right\}\eqno(6.2.2)$$
while its complement in $so(3)$ is$$
\underline{p}=\left\{\begin{pmatrix}
0&0&b\\
0&0&c\\
-b&-c&0\end{pmatrix};\ b,c\in\Bbb R\right\}.\eqno(6.2.3)$$
For the extended frame ${\cal U}^\lambda:M\to\Lambda{\rm SO}(3)_P$, with$$
{\cal U}^\lambda(0,0,\lambda)={\cal U}(0,0)=I,\eqno(6.2.4)$$

we have the Lax system ((3.2.6), restated).$$
\begin{cases}
({\cal U}^\lambda)^{-1}\cdot({\cal U}^\lambda)_x=\dis\begin{pmatrix}
0&-\varphi_x&0\\
\varphi_x&0&\lambda\\
0&-\lambda&0\end{pmatrix}={\cal A},\\
\\
({\cal U}^\lambda)^{-1}\cdot({\cal U}^\lambda)_y=\dis\begin{pmatrix}
0&0&-\lambda^{-1}\sin\varphi\\
0&0&-\lambda^{-1}\cos\varphi\\
\lambda^{-1}\sin\varphi&\lambda^{-1}\cos\varphi&0\end{pmatrix}={\cal B}.
\end{cases}\eqno(6.2.5)$$
Consequently, the Maurer-Cartan form is written as
\begin{align*}
\omega^\lambda&=-({\cal U}^\lambda)^{-1}\cdot\d{\cal U}^\lambda=-{\cal
A}\cdot\d
x-{\cal B}\cdot\d y\\
&=\alpha_{-1}\cdot\lambda^{-1}+\underbrace{(\alpha'_0\d x+\alpha''_0\d
y)}_{\alpha_0}+
\alpha_{1}\cdot\lambda,
\end{align*}
where, obviously,
\begin{align}
\alpha_{-1}&=\begin{pmatrix}
0&0&\sin\varphi\\
0&0&\cos\varphi\\
-\sin\varphi&-\cos\varphi&0\end{pmatrix}\d y,\tag{6.2.6.a}\\
\alpha_1&=\begin{pmatrix}
0&0&0\\
0&0&-1\\
0&1&0
\end{pmatrix}\d x,\tag{6.2.6.b}\\
\alpha_0'&=\begin{pmatrix}
0&\varphi_x&0\\
-\varphi_x&0&0\\
0&0&0\end{pmatrix},\qquad\alpha''_0=0.
\tag{6.2.6.c}
\end{align}

\begin{Def}
For any real smooth function $f(x,y)$ defined on a sufficiently small
neighborhood of
$(0,0)$ in $D$, we shall call $f(x,0)$ the {\it x-part} (of $f$),
respectively $f(0,y)$
the {\it y-part}.

We also set
\begin{align}
f^x&:= f(x,0)\nonumber\\
f^y&:= f(0,y)\tag{6.2.7}
\end{align}

We call $f(x,0) dx$ the {\it x-part} of the form $f(x,y) dx$.
Analogously, we call $f(0,y) dy$ the {\it y-part} of the form $f(x,y) dy$.

\end{Def}

Let $N:D\to S^2$ be the Gauss map of a weakly regular pseudospherical
surface $M$.
Thus, $N$ is real and smooth, and Lorentz harmonic. By Remark 5.1.4,
there is a $\lambda$- family of frames
${\cal U}^{\lambda }:D\to\Lambda{\rm SO}(3)_P$ such that
$\pi\circ{\cal U}^\lambda|_{\lambda=1}=N$
and such that $-({\cal U}^\lambda)^{-1}\d{\cal U}^\lambda$ is the
corresponding
Maurer-Cartan form $\omega^{\lambda}$.

Consequently, the 1-forms $\alpha_0$ and $\alpha_1$ defined by
$\omega^{\lambda}=\alpha_{-1}\lambda^{-1}+\alpha_0+\alpha_1\lambda$
are also smooth in $x$ and $y$.

In $\omega^\lambda=-({\cal U}^\lambda)^{-1}\d{\cal U}^\lambda=-{\cal A}\d
x-{\cal
B}\d y$, where ${\cal A}=({\cal U}^\lambda)^{-1}{\cal U}^\lambda_x$ and ${\cal
B}=({\cal U}^\lambda)^{-1}{\cal U}^\lambda_y$, we denote by
\begin{align}
\beta&:=\text{the $x$-part of the form $-{\cal A}\d x$ at
$\lambda=1$}\nonumber\\
       &\qquad\text{= the $x$-part of the form $\alpha'_0\d x+
\alpha_1$.}\tag{6.2.8a}\\
\beta_0&:=\text{the $x$-part of the form $\alpha'_0\,\d x$, where}\nonumber\\
       &\qquad\alpha_0(x,y)=\alpha'_0(x,y)\d x+\alpha''_0(x,y)\d
y.\tag{6.2.8b}\\
\beta_1&:=\beta-\beta_0=\text{the $x$-part of $\alpha_1$ }.\tag{6.2.8c}
\end{align}

$\gamma_0$ and $\gamma_1$ above are the analogs of $\beta_0$ and $\beta_1$,
with respect to $y$. That is

\begin{align}
\gamma&:=\text{the $y$-part of the form $-{\cal B}\d y$ at
$\lambda=1$}\nonumber\\
       &\qquad\text{= the $y$-part of $\alpha''_0 \d y+
\alpha_{-1}$.}\tag{6.2.9a}\\
\gamma_0&:=\text{the $y$-part of the form $\alpha''_0\,\d y$, where}\nonumber\\
       &\qquad\alpha_0(x,y)=\alpha'_0(x,y)\d x+\alpha''_0(x,y)\d
y.\tag{6.2.9b}\\
\gamma_1&:=\gamma-\gamma_0=\text{the $y$-part of $\alpha_{-1}$ }.\tag{6.2.9c}
\end{align}

Remark that formulas like the ones above may be also useful in the study of
other types of surfaces parametrized in some real coordinates $x,y$.
That is why we shall state some results (like Theorems 6.2.3 and 6.2.4) using
generic $\beta$'s and $\gamma$'s.

For our purposes, it is important to make $\beta$ and $\gamma$ explicit for
pseudospherical
surfaces in Chebyshev parametrization. We obtain:
$$
\beta=\beta_{1}+\beta_{0}\eqno(6.2.10a)$$
$$
\beta_1=\alpha_1=\begin{pmatrix}
0&0&0\\
0&0&-1\\
0&1&0\end{pmatrix}\d x,\eqno(6.2.10b)$$
$$
\beta_0=\alpha_0'(x,0)\d=\begin{pmatrix}
0&{\varphi_x}(x,0)&0\\
{-\varphi_x}(x,0)&0&0\\
0&0&0\end{pmatrix}\d x.\eqno(6.2.10c)$$
$$
\gamma=\gamma_{1},\eqno(6.2.11a)$$
$$
\gamma_{1}=\begin{pmatrix}
0&0&\sin\varphi(0,y)\\
0&0&\cos\varphi(0,y)\\
-\sin\varphi(0,y)&-\cos\varphi(0,y)&0\end{pmatrix}\d y,\eqno(6.2.11b)$$
$$
\gamma_0=0.\eqno(6.2.11c)$$

Let us now recall the two Birkhoff-type factorizations presented in Theorem
4.3.2.

The first type of Birkhoff factorization from Theorem 4.3.2 is performed on
the ``big cell" $\tilde\Lambda_*^-{\rm SO}(3)_P\cdot\tilde\Lambda^+{\rm
SO}(3)_P$.
That is, away from a singular set $S_1\subset D$, we can split the extended
moving frame
${\cal U}^\lambda:D\to{\rm SO}(3)$ into two parts.
Recall that the first factor of this splitting is of the form
$g_-=I+\lambda^{-1}g_{-1}+\lambda^{-2}g_{-2}+\cdots$, while the second
factor of
the splitting is of the form $g_+=g_0+\lambda g_1+\lambda^2g_2+\cdots$,
respectively.

Since the ``big cell" is open and ${\cal U}^\lambda:D\to{\rm SO}(3)$ is
continuous,
the set $$
{\tilde D_1}=\{(x,y)\text{ ; }{\cal U}^{\lambda}(x,y)\text{belongs to the
``big cell"}\}$$
is open.
Note that $(0,0)\in\tilde D_1$.

Let $S_1 = D -\tilde D_1$ denote the ``singular" set. We have just shown
that $S_1$
is closed and $(0,0)$ is not an element of the set $S_1$. Similarly, we
have $S_2$ and
$\tilde D_2$ for the second splitting.

The second type of Birkhoff splitting is the analogous splitting in the
``big cell" $\tilde\Lambda_*^+{\rm SO}(3)_P\times\tilde\Lambda^-{\rm SO}(3)_P$.
The goal of this section is to show that the first factor of each type of
splitting is
an essential one, and can be viewed as an integral of the unconstrained
data that we call
normalized potential.

We can perform the two splittings on the extended frame $\cal U^\lambda$.
{\it Let ${\cal U}={\cal U}^\lambda$ be the
extended normalized moving frame of a pseudospherical surface and let $(x,y)\in
D\setminus(S_1\cup S_2)$. Then, for some uniquely determined
$V_+\in\Lambda^+{\rm
SO}(3)_P$, $V_-\in\Lambda^-{\rm SO}(3)_P$ and
${\cal U}_-\in\Lambda_*^-{\rm SO}(3)_P$, ${\cal U}_+\in\Lambda_*^+{\rm
SO}(3)_P$,
${\cal U}$ can be written as}$$
{\cal U}={\cal U}_+\cdot V_-={\cal U}_-\cdot V_+.\eqno(6.2.12)$$
\vskip6pt

The factors carrying the ``genetic material" to recreate the frame and then the
surface are ${\cal U}_+$ and ${\cal U}_-$. They can be obtained, starting
from two
normalized potentials $\eta ^x$ and $\eta ^y$ respectively, by
solving the two ordinary differential equations presented in Theorem 6.2.1.

From $\cal U_-$ and $\cal U_+$, one can reproduce the frame $\cal U$ and
then construct the corresponding pseudospherical immersion via the
Sym-Bobenko formula
(5.1.5).

\begin{Th}
Let ${\cal U}={\cal U}^\lambda$, ${\cal U}_+$ and ${\cal U}_-$ be
as above. Then the following systems of differential equations are
satisfied:
 {\rm1)}\hfill$({\cal U}_+)^{-1}\dis{\partial{\cal
U}_+\over\partial x} \d x=-\lambda\cdot V_0\cdot\beta_{1}\cdot
V_0^{-1}$,\hfill$(6.2.13)$ \vskip6pt \noindent with initial
condition ${\cal U}_+(x=0)=I$, where $V_0$ is some matrix
$V_0(x)\in{\rm SO}(3)$. \vskip6pt {\rm2)}\hfill$({\cal
U}_-)^{-1}\dis{\partial{\cal U}_-\over\partial y} \d
y=-\lambda^{-1}\cdot W_0\cdot\gamma_{1}\cdot
W_0^{-1}$,\hfill$(6.2.14)$ \vskip6pt \noindent with initial
condition ${\cal U}_-(y=0)=I$, where $W_0$ is some matrix
$W_0(y)\in{\rm SO}(3)$.

Moreover, ${\cal U}_+$ does not depend on $y$ and ${\cal U}_-$ does not
depend on
$x$.

In some other words, ${\cal U}_+$ and ${\cal U}_-$ are solutions of some
first order systems of differential equations in $x$ and $y$, respectively.
\end{Th}

\noindent{\it Proof. of 2)} From equation (6.2.12), we know$$
{\cal U}_-={\cal U}\cdot V_+^{-1}.\eqno(6.2.15)$$
Differentiating (6.2.15), we obtain$$
\d{\cal U}_-=\d{\cal U}\cdot V_+^{-1}-{\cal U}\cdot V_+^{-1}\cdot\d V_+\cdot
V_+^{-1},\eqno(6.2.16)$$
which can be rewritten as$$
{\cal U}_-^{-1}\cdot d{\cal U}_-=V_+\cdot({\cal U}^{-1}\cdot\d{\cal U})\cdot
V_+^{-1}-\d V_+\cdot V_+^{-1}\eqno(6.2.17)$$
after left multiplication by ${\cal U}_-^{-1}$.

The coefficient of $\d x$ on the left-hand side of (6.2.17) contains only
negative powers of
$\lambda$, while the coefficient of $\d x$ on the right-hand side of (6.2.17),
in view of
(6.2.5), contains only non-negative powers of $\lambda$. Therefore,
$\partial_x{\cal U}_-=0$, so
${\cal U}_-$ depends on $y$ only.

To determine (6.2.14), we consider the coefficient of $\d y$ in (6.2.17). The
left-hand side of
(6.2.17) contains only negative powers of $\lambda$, while the one on the
right-hand side,
due to (6.2.5),$$
{\cal B}={\cal U}^{-1}\partial_y{\cal U}=\lambda^{-1}\cdot\begin{pmatrix}
0&0&-\sin\varphi\\
0&0&-\cos\varphi\\
\sin\varphi&\cos\varphi&0
\end{pmatrix},$$
contains only one term in $\lambda^{-1}$, and no terms in $\lambda^k$, $k<-1$.

On the other hand, let$$
V_+=\tilde{W}_0+\lambda\tilde{W}_1+\lambda^2\tilde{W}_2+\cdots=\tilde{W}_0\cdot
T_+,\eqno(6.2.18)$$
with $T_+\in\Lambda^+_*{\rm SO}(3)_P$.

Therefore, ${\cal U}_-^{-1}\partial_y{\cal U}_-=\tilde{W}_0{\cal
B}\tilde{W}_0^{-1}$, where the left-hand side only depends on $y$. Since ${\cal
U}$,
$V_+$ and $\tilde{W}_0$ are all defined on
$D$, a neighborhood of $(0,0)$, we can specialize to the points of the form
$(0,y)$
for a sufficiently small interval on the line
$x=0$, containing the origin.

Thus,$$
{\cal U}_-^{-1}\partial_y{\cal U}_-=\tilde{W}_0(0,y)\cdot{\cal
B}(0,y)\cdot\tilde{W}_0(0,y)^{-1}.\eqno(6.2.19)$$

We observe that$$
{\cal B}(0,y)=\lambda^{-1}\cdot\begin{pmatrix}
0&0&-\sin{\varphi (0,y)}\\
0&0&-\cos{\varphi (0,y)}\\
\sin{\varphi (0,y)}&\cos{\varphi (0,y)}&0
\end{pmatrix}$$

From formulas (6.2.5-6), we note that ${\cal B}=-\lambda^{-1}\cdot
\alpha_{-1}$,
and restricting to the $y$-parts, we obtain  $$
{\cal B}(0,y) \d y=-\lambda^{-1}\cdot {\gamma_1},$$
where the form $\gamma_1$ is the one given in formulas (6.2.11.c).

In (5.2.2), we defined the normalized $y$-potential by
$$
\eta^y= -{\lambda}\cdot {\cal U}_-^{-1}\partial_y{\cal U}_- \d
y.\eqno(6.2.20)$$
and the meromorphic $y$-potential as$$
\xi^y= - {\cal U}_-^{-1}\partial_y{\cal U}_- \d y=\lambda^{-1}\eta^y,$$

Denoting $W_0(y):=\tilde{W}_0(0,y)$, we obtain$$
{\cal U}_-^{-1}\partial_y{\cal U}_- \d y= -\lambda^{-1}\cdot
W_0(y)\cdot\gamma_{1}\cdot[W_0(y)]^{-1}$$

and therefore (6.2.14).
\vskip6pt
\noindent{\it Proof of 1)}

From equation (6.2.12), we obtain$$
{\cal U}_+={\cal U}\cdot V_-^{-1},\quad {\cal U}_+\in\Lambda_*^+{\rm
SO}(3)_P,\quad
V_-\in\Lambda^-{\rm SO}(3)_P,\eqno(6.2.21)$$
which by differentiation leads to$$
\d{\cal U}_+=\d{\cal U}\cdot V_-^{-1}-{\cal U}\cdot V_-^{-1}\cdot\d V_-\cdot
V_-^{-1},\eqno(6.2.22)$$
and then
\vskip-24pt$$
{\cal U}_+^{-1}\d{\cal U}_+=V_-({\cal U}^{-1}\d{\cal U})V_-^{-1}-\d V_-\cdot
V_-^{-1}.\eqno(6.2.23)$$
We compare the coefficient of $\d y$ on the left-hand side of (6.2.23) with
the coefficient
of $\d y$ on the right-hand side of (6.2.23), via formula (6.2.5). The
left-hand side
of (6.2.23) clearly contains only positive powers of $\lambda$, while the
coefficient of
$\d y$ on the right-hand side of (6.2.23), in view of (6.2.5), contains
non-positive powers
of $\lambda $ only.
Thus, ${\cal U}_+$ depends exclusively on $x$.

In order to obtain (6.2.13), we consider the coefficient of $\d x$ in
(6.2.23).
The left-hand side of (6.2.23) contains only positive powers of $\lambda$,
while the one
on the right-hand side, due to
$$
{\cal A}={\cal U}^{-1}\partial_x{\cal U}=\begin{pmatrix}
0&-\varphi_x&0\\
\varphi_x&0&\lambda\\
0&-\lambda&0\end{pmatrix},$$
contains one term in $\lambda$ and no terms in $\lambda^k$, with $k>1$.

Like we did before in case 2), we can restrict to a sufficiently small interval
around
$(0,0)$ on the line $y=0$.

Let now$$
V_-=\tilde{V}_0+\lambda^{-1}\tilde{V}_1+\lambda^{-2}\tilde{V}_2+\cdots=\tilde{V}
_0\cdot
T_-,\eqno(6.2.24)$$
with $T_-\in\Lambda_*^-{\rm SO}(3)_P$.

Then we note$$
{\cal U}_+^{-1}(x)\cdot\partial_x{\cal U}_+=\tilde{V}_0(x,0)\cdot{\cal
A}(x,0)\cdot\tilde{V}_0(x,0)^{-1}.\eqno(6.2.25)$$
Moreover, since the left-hand side of (6.2.23) contains only positive
powers of $\lambda$,
we conclude that$$
{\cal U}_+^{-1}(x)\cdot\partial_x{\cal U}_+ \d
x=-\tilde{V}_0(x,0)\cdot\lambda\cdot
\beta_{1}\cdot\tilde{V}_0(x,0)^{-1},\eqno(6.2.26)$$
where according to formula (6.2.10.b), $\beta_{1}=\alpha_{1}= - E_{23}$.
\
This is exactly the claim of the equation (6.2.13) stated in the theorem,
if we denote $\tilde{V}_0(x,0) := V_0 $.\hfill$\square$

\subsection{Normalized Potentials for Pseudospherical Surfaces}

In this section we find the explicit expressions of the two normalized
potentials.
Theorems 6.3.1. and 6.3.2 can be thought of as corollaries to Theorem 6.2.1.
Basically, we construct the normalized potentials from
the solutions to the ordinary differential systems introduced in Theorem 6.2.1.
Theorems 6.3.1. and 6.3.2 are phrased analogously to Wu's Theorem 6.1.1 for
the normalized
potential of the constant mean curvature surfaces.

\begin{Th} {} ($x$-potential) The normalized potential $\eta^x$  with the
origin as the reference point is given by$$
\eta^x=V_0(x)\cdot\beta_1(x)\cdot V_0(x)^{-1},\eqno(6.3.1)$$
where $V_0$ is the solution of$$
\begin{cases}
V_0(x)^{-1}\d V_0(x)=-\beta_0(x),\\
V_0(0)={\cal U}(0,0).
\end{cases}
\eqno(6.3.2)$$ where $\beta_0$ and $\beta_1$ are given by formulas
(6.2.10). \vskip-6pt

\end{Th}

Similarly, we have the following result:

\begin{Th} ($y$-potential) The normalized potential $\eta^y$ with the
origin as the
reference point is given by$$
\eta^y=W_0(y)\cdot\gamma_1(y)\cdot W_0(y)^{-1},\eqno(6.3.3)$$
where $W_0$ is the solution of$$
\begin{cases}
W_0(y)^{-1}\d W_0(y)=-\gamma_0(y),\\
W_0(0)={\cal U}(0,0).
\end{cases}
\eqno(6.3.4)$$
where $\gamma_0$ and $\gamma_1$ are given by formulas (6.2.11).
\vskip-6pt

\end{Th}

\

{\it Proof of Theorem 6.3.1}.

Relation (6.3.3) is a rephrasing of (6.2.13):$$ {\cal
U}_+^{-1}\partial_x{\cal U}_+ \d x= -\lambda\cdot
V_0(x)\cdot\beta_{1}\cdot V_0(x)^{-1},$$ where we substitute
$\eta^x = -\lambda^{-1}\cdot ({\cal U_+})^{-1}\cdot \d {\cal
U_+}$, that is the definition of the normalized $x$-potential.

Let us now consider again equation (6.2.23) from the proof of Theorem
6.2.1, namely$$
{\cal U}_+^{-1}\d{\cal U}_+=V_-({\cal U}^{-1}\d{\cal U})V_-^{-1}-\d V_-\cdot
V_-^{-1}$$
We proved that both sides depend on $x$ only. Now let us take a look at the
coefficient of $\lambda^{0}$ in this equation.

The left-hand side has positive powers of $\lambda$ only, while the
$x$-part of
right-hand side only has $-V_0\cdot\beta_{0}\cdot V_0^{-1}-\d V_0\cdot
{V_0}^{-1}$
as a term that does not depend on $\lambda$.

Consequently, we obtain $V_0(x)^{-1}\d V_0=-\beta_0(x)$. Formula (6.2.10.c)
shows
that$$
\beta_0=\alpha_0'(x,0)=\begin{pmatrix}
0&{\varphi_x}(x,0)&0\\
-{\varphi_x}(x,0)&0&0\\
0&0&0\end{pmatrix}.$$
Here it was taken into account that $\varphi_x(x,0)=(\varphi(x,0))_x$, where
$\xi(x):=\varphi(x,0)$ is the part in $x$ of the smooth angle function
$\varphi(x,y)$.
If we consider the matrix$$
\theta=\begin{pmatrix}
0&\xi&0\\
-\xi&0&0\\
0&0&0\end{pmatrix},\eqno(6.3.5)$$
then $$
\beta_0=\theta'\d x=\begin{pmatrix}
0&\xi'(x)&0\\
-\xi'(x)&0&0\\
0&0&0\end{pmatrix}\d x.\eqno(6.3.6)$$
The solution $V_0$ of the system (6.3.2) must take into account that
${\cal U}(0,0,\lambda)=I$, so the solution is$$
V_0(x)=e^{\theta(0)-\theta(x)}.\eqno(6.3.7)$$
Using also the expression of the form
$\beta_1$, the normalized $x$-potential $\eta^x$ can be written as$$
\eta^x=V_0(x)\beta_1(x)V_0^{-1}(x)=e^{\theta(0)-\theta(x)}(-E_{23})e^{\theta(x)-
\theta(0)}\d x.\eqno$$
Since$$
e^{\theta(0)-\theta(x)}=\begin{pmatrix}
\cos(\xi(0)-\xi(x))&\sin(\xi(0)-\xi(x))&0\\
-\sin(\xi(0)-\xi(x))&\cos(\xi(0)-\xi(x))&0\\
0&0&1
\end{pmatrix},\eqno(6.3.8)$$
the formula above leads to the final expression of the $x$-potential, as$$
\eta^x=\begin{pmatrix}
0&0&-\sin(\xi(0)-\xi(x))\\
0&0&-\cos(\xi(0)-\xi(x))\\
\sin(\xi(0)-\xi(x))&\cos(\xi(0)-\xi(x))&0
\end{pmatrix}\d x,\eqno(6.3.9)$$
where $\xi(x):=\varphi (x,0)$.
\mbox{ }\hfill$\square$

{\it Proof of Theorem 6.3.2}.

Relation (6.3.3) is a rephrasing of (6.2.14):$$
{\cal U}_-^{-1}(y)\cdot\partial_y{\cal U}_- =-\lambda^{-1}\cdot
W_0\cdot\gamma_{1}\cdot W_0^{-1}$$
where we substitute (6.2.20)$$
\eta^y= -{\lambda}\cdot {\cal U}_-^{-1}\partial_y{\cal U}_-,$$
that is the definition of the normalized $y$-potential.

Let us now consider equation (6.2.17) from the proof of Theorem 6.2.1,
namely$$
{\cal U}_-^{-1}\cdot d{\cal U}_-=V_+\cdot({\cal U}^{-1}\cdot\d{\cal U})\cdot
V_+^{-1}-\d V_+\cdot V_+^{-1}.$$

We proved that both sides depend on $y$ only. Now let us take a look at the
coefficient of $\lambda^{0}$ in this equation.

The left-hand side has negative powers of $\lambda$ only, while the
$y$-part of
right-hand side only has $-W_0\cdot\gamma_{0}\cdot W_0^{-1}-\d W_0\cdot
{W_0}^{-1}$
as a term that does not depend on $\lambda$.

Consequently, we obtain $W_0(y)^{-1}\d W_0=-\gamma_0(x)$. Formula (6.2.11.c)
tells us that $\gamma_0=0$.
From this we conclude that $W_0(y)$ is actually a constant matrix, and from
the initial condition on the frame $\cal U$, together with the initial
condition of (6.3.4), it follows that for every $y$,$$
W_0(y)={\cal U}(0,0)=I.$$
It follows that$$
\eta^y=W_0(y)\cdot\gamma_1(y)\cdot W_0(y)^{-1}=\gamma_1(y).$$

Therefore,$$
\eta^y=\begin{pmatrix}
0&0&\sin\varphi(0,y)\\
0&0&\cos\varphi(0,y)\\
-\sin\varphi(0,y)&-\cos\varphi(0,y)&0\end{pmatrix}\d y.\eqno(6.3.10)$$

\mbox{ }\hfill$\square$

\vskip12pt
\begin{Rem} Let us review the
expressions (6.3.9) and (6.3.10) for the two normalized potentials $\eta^x$
and $\eta^y$, that is$$
\begin{pmatrix}
0&0&-\sin(\varphi(0,0)-\varphi(x,0))\\
0&0&-\cos(\varphi(0,0)-\varphi(x,0))\\
\sin(\varphi(0,0)-\varphi(x,0))&\cos(\varphi(0,0)-\varphi(x,0))&0
\end{pmatrix}\d x,$$
and$$
\begin{pmatrix}
0&0&\sin(\varphi(0,y))\\
0&0&\cos(\varphi(0,y))\\
-\sin(\varphi(0,y))&-\cos(\varphi(0,y))&0
\end{pmatrix}\d y,$$
respectively.

Note that the normalized potentials depend exclusively on the angle
$\varphi(x,y)$ between the asymptotic lines.

\end{Rem}

\subsection{Another Method for Normalized Potentials: Passage to 2 $\times$ 2
Matrices}

We introduce the matrices $$
\sigma_1=\begin{pmatrix}
0&1\\1&0\end{pmatrix},\quad
\sigma_2=\begin{pmatrix}
0&-i\\i&0\end{pmatrix},\quad
\sigma_3=\begin{pmatrix}
1&0\\0&-1\end{pmatrix}\eqno(6.4.1)$$
called {\it Pauli matrices}.

We can rewrite in terms of $2\times2$ matrices the potentials $\eta^x$ and
$\eta^y$, given by
formulas (6.3.9) and (6.3.10). The Pauli matrices above will allow us to do
this.

This passage from $3\times3$ to $2\times2$ matrices can be done via the
well-known isomorphism between
${\rm sl}(2,\Bbb C)$ and ${\rm so}(3,\Bbb C)$, which induces an isomorphism
between ${\rm su}(2)$ and ${\rm so}(3,\Bbb R)$. This isomorphism is defined by
the correspondence
$$
E_{12}\longleftrightarrow (-i/2)\sigma_3,\eqno(6.4.2a)$$$$
E_{13}\longleftrightarrow (-i/2)\sigma_2,\eqno(6.4.2b)$$$$
E_{23}\longleftrightarrow (-i/2)\sigma_1.\eqno(6.4.2c)$$

Via this passage to $2\times2$ matrices, the two potentials become
\begin{align}
\eta^x&=\frac{i}{2}\begin{pmatrix}
0&e^{i(\varphi(x,0)-\varphi(0,0))}\\
e^{i(\varphi(0,0)-\varphi(x,0))}&0\end{pmatrix}\d x\tag{6.4.3a},\\
\eta^y&=-\frac{i}{2}\begin{pmatrix}
0&e^{-i(\varphi(0,y))}\\
e^{i(\varphi(0,y))}&0\end{pmatrix}\d y,\tag{6.4.3b}
\end{align}
respectively.

\begin{Rem}
The fact that we chose to work with Chebyshev nets ($A=B=1$, where
$A=|\psi_x|$, $B=|\psi_y|$) allows this form of the normalized pair of
potentials.

Had we chosen to work with arbitrary $A$ and $B$ (that is, not necessarily a
Chebyshev net), a
more tedious but straightforward calculation would lead us to the
Weierstrass data
\begin{align*}
\eta^x&=\frac{iA}{2}\begin{pmatrix}
0&e^{i(\varphi(x,0)-\varphi(0,0))}\\
e^{i(\varphi(0,0)-\varphi(x,0))}&0
\end{pmatrix}\d x,\tag{6.4.4a}\\
\eta^y&=-\frac{iB}{2}\begin{pmatrix}
0&e^{-i(\varphi(0,y))}\\
e^{i(\varphi(0,y))}&0
\end{pmatrix}\d y.\tag{6.4.4b}
\end{align*}
\end{Rem}

The asymmetry of the two potential comes from the definition of a
normalized extended frame. Although the two potentials look
asymmetric, one can make the two expressions look similar by
gauging with a certain rotation.

\begin{Rem} The product of the off-diagonal elements is $A^2$ and
$B^2$ respectively for $\eta^x$ and $\eta^y$ (with a factor of $-1/4$).
This is similar to the CMC case, where the meromorphic (normalized)
potential has the form$$
\eta=\begin{pmatrix}
0&f(z)\\
g(z)&0\end{pmatrix}
\d z,\eqno(6.4.5)$$
with $f\cdot g \d z^2=-Q \d z^2$ (Hopf differential).

For the CMC case, the $\lambda$-transformation was given by $$
\begin{cases}
Q\mapsto e^{2it}Q=\lambda^2Q,\\
\bar{Q}\mapsto e^{-2it}\bar{Q}=\lambda^{-2}Q,
\end{cases}\eqno(6.4.6)$$
while here it is $A\mapsto\lambda A$, $B\mapsto\lambda^{-1}B$,
$\lambda=e^t$. So
the role played in the case of CMC surfaces by the Hopf differential $Q$ is
taken
for the case of pseudospherical surfaces by the pair $A,B$.

The globally defined differential forms $(A^2) \d x^2$ and $(B^2) \d y^2$
are sometimes called Klotz differentials.

\end{Rem}

\begin{Rem}
The isomorphism described above in (6.4.2 a,b,c), between ${\rm su}(2)$ and
${\rm so}(3)$,
is provided by the spinor representation J defined as follows:$$
J:\Bbb R^3\to{\rm su}(2),$$$$
J(x,y,z)=\frac{1}{2}\begin{pmatrix}
-iz&-ix-y\\
-ix+y&iz\end{pmatrix},\eqno(6.4.7)$$
which identifies $\Bbb R^3$ and ${\rm su}(2)$ via$$
J({\bf r})=-\frac{i}{2}{\bf r}\sigma,\eqno(6.4.8)$$
where ${\bf r}\sigma=r_1\sigma_1+r_2\sigma_2+r_3\sigma_3$, and
$\sigma_1,\sigma_2,\sigma_3 $
are the Pauli matrices defined by (6.4.1). Then$$
J({\bf r}_1\times{\bf r}_2)=[J{\bf r}_1,J{\bf r}_2].\eqno(6.4.9)$$

\end{Rem}

If ${\cal U}=(e_1,e_2,e_3)$ is the normalized moving frame of the surface
$M$ in
asymptotic line parametrization, we define the ``$2\times2$" {\it frame}
$P:D\to{\rm SU}(2)$ with the initial condition $P(0,0)=I$, via
\begin{align}
J(e_1)&=-\frac{i}{2}P\sigma_1P^{-1},\nonumber\\
J(e_2)&=-\frac{i}{2}P\sigma_2P^{-1},\nonumber\\
J(e_3)&=-\frac{i}{2}P\sigma_3P^{-1}.\tag{6.4.10}
\end{align}
  We have this way a correspondence between all the frames ${\cal U}$ in
${\rm SO}(3)$ and frames $P$ in ${\rm SU}(2)$.

A tedious but straightforward computation completely similar to the one in
[DoHa], Appendix A.4, transfers the $3\times3$ matrices ${\cal A}$ and
${\cal B}$
from
(6.2.5) to the $2\times 2$ matrices $U$ and $V$ given by the corresponding
Lax
system:
$$
U=P^{-1}P_x=\frac{-i}{2}\begin{pmatrix}
-\varphi_x&\lambda\\
\lambda &\varphi_x\end{pmatrix}\eqno(6.4.11)$$
$$
V=P^{-1}P_y=\frac{i}{2}\lambda^{-1}\begin{pmatrix}
0&e^{-i\varphi}\\
e^{i\varphi}&0\end{pmatrix}\eqno(6.4.12)$$
with $P(0,0)=I$.

This Lax system can be also obtained directly from (2.3.15) through the
isomorphism
between ${\rm su}(2)$ and ${\rm so}(3,\Bbb R)$ defined by (6.4.2 a,b,c),
that is
$$
E_{12}\longleftrightarrow (-i/2)\sigma_3,\quad
E_{13}\longleftrightarrow (-i/2)\sigma_2,\quad
E_{23}\longleftrightarrow (-i/2)\sigma_1.$$

\subsection*{Appendix}

By definition, the deformation parameter $\lambda$ that
generates an associated family of pseudospherical surfaces is
real and positive, $\lambda=e^t$. In general, $\lambda$ can be
real and negative as well. Our choice is motivated by the
convenience of working within the connected Banach loop group$$
(\Lambda{\rm SO}(3)_P,\|\cdot\|)=\{g:\Bbb R_+\to{\rm SO}(3,\Bbb
R)\mid Pg(\lambda)P^{-1}=g(-\lambda)\},$$
endowed with the norm $\|\cdot\|$ defined by (4.1.16).

The goal of this appendix is to show that we can split \`a  la
Birkhoff any extended frame ${\cal U}^\lambda$ which admits an analytic
extension on $\Bbb C_*$.

For our purpose, it is useful to extend the real positive parameter
$\lambda$ such that we can apply the Birkhoff splitting to complex loop groups
with loop parameter in $S^1$.

Let us first consider the Lax system$$
\begin{cases}
{\cal U}^{-1}\cdot\partial_x{\cal U}=-\varphi_x\cdot E_{12}+\mu\cdot E_{23}\\
{\cal U}^{-1}\cdot\partial_y{\cal U}=\mu^{-1}\cdot(-\sin\varphi\cdot
E_{13}-\cos\varphi\cdot E_{23}),
\end{cases}\eqno({\rm A}.1)$$
where $\mu\in\Bbb C_*=\Bbb C\setminus\{0\}$.

Clearly, the Lax system (2.3.15) is the same as ({\rm A}.1) if we
restrict $\mu$ to $\Bbb R_+$. \vskip6pt \noindent{\bf Lemma A.1}
\it Every solution $\cal U$ to {\rm(A.1)} with initial condition
${\cal U}(0,0,\mu)=I$ is analytic in $\mu\in\Bbb C_*$ and$$
\overline{{\cal U}(x,y,\bar{\mu})}={\cal U}(x,y,\mu).\eqno({\rm
A}.2)$$ \vskip6pt \noindent Proof. \rm Note that the right-hand
side of $A.1$ is analytic in $\mu\in\Bbb C_*$. Since the initial
condition is analytic in $\mu\in\Bbb C_*$, it follows that, for
every $x$ and $y$ arbitrarily fixed, the solution $\cal U$ of
(A.1) is also analytic in $\mu\in\Bbb C_*$. Relation (A.2) is
straight-forward, as a consequence of the reality of
$\varphi$.\hfill$\square$

In order to use the classical loop group factorization, let us
choose $\lambda\in S^1$. We consider the restriction to $S^1$ of the
extended frame $\cal U$ satisfying the Lax system (A.1), and will
denote it $\cal U^\lambda$.

Taking into consideration the property (A.2) of $\cal U^\lambda$,
we introduce the
following group of continuous maps:
\begin{align}
H_P=\{A:S^1\to{\rm SO}(3,\Bbb C)\mid\ &A\text{ continuous,
}\overline{A(\bar{\lambda})}=A(\lambda),\nonumber\\
&P\cdot A(\lambda)\cdot P^{-1}=A(-\lambda)\},
\tag{{\rm A}.3}
\end{align}
with the supplementary condition$$
\|A\|=\sum_{k\in\Bbb Z}\|A_k\|<\infty,\eqno({\rm A}.4)$$
where$$
A(\lambda)=\sum_{k\in\Bbb Z}A_k\cdot\lambda^k,$$
and$$
\|B\|=\max_{i} \{ {\sum_{j=1}^3|B_{ij}|} \},\eqno({\rm A}.5)$$

for every $\lambda$-independent $3\times 3$ matrix $B$.

$H_P$ is a Banach Lie group with respect to the norm $\|\cdot\|$.

(Note that we have used the same symbol $\|\cdot\|$ for different entities.
We hope this will not lead to any confusion).

Clearly, (A.3) expresses the reality of the coefficient matrices in
the Fourier expansion.

\noindent{\bf Proposition A.1} {\it
 For the group $H_P$ we define $(H^-_P)_*$ and $H_P^+$ as in Section 4.3.
The multiplication $(H^-_P)_*\times H_P^+\to (H^-_P)_*\cdot H^+_P$
is an analytic diffeomorphism onto the open and dense subset
$(H^-_P)_*\cdot H^+_P$, called the ``big cell". In particular, if
$A\in H_P$ is contained in the big cell, then $A$ has a unique
decomposition$$ A=A_-A_+,\eqno({\rm A}.6)$$ where $A_-\in
(H^-_P)_*$ and $A_+\in H_P^+$. The analogous result holds for the
multiplication map $(H^+_P)_*\times H^-_P\to (H^+_P)_*\cdot
H^-_P$.}

\noindent Proof. \rm Let $A:S^1\to{\rm SO}(3,\Bbb C)$ be an
element of $H_P$. By the definition of $H_P$, we have$$
A(\lambda)=\overline{(A(\bar{\lambda})},\qquad\text{for every
}\lambda\in S^1.$$ On the other hand, by Theorem 4.3.1,
$A(\lambda)$ can be decomposed \`a la Birkhoff in a big cell of
$\Lambda{\rm SO}(3,\Bbb C)$, as $$
A(\lambda)=A_-(\lambda)A_+(\lambda),\qquad\text{for every
}\lambda\in S^1,\eqno({\rm A}.7),$$
$A_-(\lambda)\in\Lambda_*^-{\rm SO}(3,\Bbb C)$,
$A_+(\lambda)\in\Lambda^+{\rm SO}(3,\Bbb C)$.

As a consequence of (A.6) and (A.7), we obtain$$
A(\lambda)=\overline{A_-(\bar{\lambda})}\cdot\overline{A_+(\bar{\lambda})}.
\eqno({\rm A}.8)$$

Also, (A.7) and (A.8) yield$$
A_-(\lambda)^{-1}\cdot\overline{A_-(\bar{\lambda})}=
A_+(\lambda)\cdot\overline{A_+(\bar{\lambda})}^{-1}.$$

The left-hand side is an element of $\Lambda_*^-{\rm SO}(3,\Bbb C)$, while the
right-hand side is an element of $\Lambda^+{\rm SO}(3,\Bbb C)$, and hence
both sides are
equal to the identity matrix. Therefore,
\begin{align}
\overline{A_-(\bar{\lambda})}&=A_-(\lambda)\nonumber\\
\overline{A_+(\bar{\lambda})}&=A_+(\lambda),\qquad\text{for every }\lambda\in
S^1.\tag{{\rm A}.9}
\end{align}

Hence, $A_+$ and $A_-$ satisfy the first condition (A.6) from the
definition of the
group $H_P$, meaning that their coefficient matrices are real.

On the other hand, the symmetry condition$$
A(-\lambda)=P\cdot A(\lambda)\cdot P^{-1}=P\cdot A_-(\lambda)\cdot P^{-1}\cdot
 P\cdot A_+(\lambda)\cdot P^{-1},$$
together with the uniqueness of the Birkhoff splitting$$
A(-\lambda)=A_-(-\lambda)\cdot A_+(-\lambda),$$
yield the symmetry condition for $A_-$, $A_+$:
\begin{align*}
A_-(-\lambda)&=P\cdot A_-(\lambda)\cdot P^{-1}\\
A_+(-\lambda)&=P\cdot A_+(\lambda)\cdot P^{-1}.
\end{align*}
Thus, the Birkhoff factorization holds for $H_P$. The analytic
diffeomorphism $(H^-_P)_*\times H_P^+\to (H^-_P)_*\cdot H^+_P$ is
a particularization of the analytic diffeomorphism analyzed in
Theorem 4.2.1. \hfill$\square$

Theorem 4.3.2. states the Birkhoff splitting for arbitrary elements
$g$ in the Banach loop group $\Lambda{\rm SO}(3)_P$ which admit an
analytic extension to $\Bbb C_*$.

Now we are ready to present its proof.

\vskip6pt
\noindent{\it Proof of Theorem 4.3.2.} Let $g:\Bbb R_+\to{\rm SO}(3)$,
$g(-\lambda)=P\cdot g(\lambda)\cdot P^{-1}$, be an element of the Banach
loop group
$\Lambda{\rm SO}(3)_p$ that has an analytic extension $\tilde{g}$ to
$\Bbb C_*$.

Set $A:=\tilde{g}|_{S^1}$.

Since $g\in\Lambda{\rm SO}(3)_P$, the matrix coefficients of $g$ are
real, that is$$
A(\lambda)=\sum_{k\in\Bbb Z}A_k\lambda^k,\qquad \lambda\in S^1,$$

Note that $A\in H_P$, where $H_P$ denotes the loop group defined by (A.3).
The algebraic conditions are obviously satisfied. Also, by [GO], Theorem 1.4,
analytic functions satisfy the finite norm condition.

Here we are only interested in elements $A$ belonging to the big cell of $H_P$.

The previous proposition shows that the Birkhoff splitting holds for the
big cell
of $H_P$.

Then, let $A_-\in (H^-_P)_*$ and $A_+\in H^+_P$ be such that$$
A=A_- A_+,$$
where
\begin{align*}
A_-&=I+\sum_{k<0}A_k\lambda^k,\\
A_+&=\sum_{k\geq0}A_k\lambda^k,\qquad\lambda\in S^1.
\end{align*}

We need to show that $A_-$ and $A_+$ admit analytic extensions to $\Bbb C_*$.

By our hypothesis, $A$ has an analytic extension to $\Bbb C_*$.
The element $A_-$ admits an analytic extension to the exterior of the unit
circle $S^1$. Therefore, $(A_-)^{-1}A = A_+$ can be extended analytically
outside of
the unit disk.

On the other hand, $A_+$ admits an analytic extension inside the unit disk.
Thus, by analytic prolongation, $A_+$ admits an analytic extension to $\Bbb
C_*$.

From  $A (A_+)^{-1}=A_- $, it follows next that $A_-$ also admits an
analytic extension
to $\Bbb C_*$.

Let $\tilde{A}_-$ and $\tilde{A}_+$ be the analytic extensions of $A_-$ and
$A_+$ to
$\Bbb C_*$, respectively.

Next, let $g_-$ and $g_+$ denote their restrictions to $\Bbb R_+$:$$
g_-=\tilde{A}_-|_{\Bbb R_+},\quad g_+=\tilde{A}_+|_{\Bbb R_+}.$$
Clearly, $g$, $g_-$ and $g_+$ have analytic extensions to $\Bbb
C_*$, respectively:
$\tilde{g}$, $\tilde{A}_-$ and $\tilde{A}_+$ such that$$
\tilde{g}|_{S^1}=A=A_-\cdot A_+=\tilde{A}_-|_{S^1}\cdot\tilde{A}_+|_{S^1},$$
that is $\tilde{g}$ and $\tilde{A}_-\tilde{A}_+$ coincide on $S^1$. Therefore,
$\tilde{g}$ and $\tilde{A}_-\tilde{A}_+$ will coincide on $\Bbb R_+$
as well, and
$g=g_-g_+$ is a unique factorization.

This proves the splitting.

It remains to prove that
$\tilde\Lambda_*^-SO(3)_P\times\tilde\Lambda^+SO(3)_P\to\tilde\Lambda SO(3)_P$
is a diffeomorphism onto the open and dense subset
$\tilde\Lambda_*^-SO(3)_P\cdot\tilde\Lambda^+SO(3)_P$.

Note that $\tilde\Lambda SO(3)_P$ is a subgroup of $\Lambda SO(3)_P$ with
the induced
topology. On the other hand, it is natural to view the diffeomorphism
$\Lambda_*^-SO(3)_P\times\Lambda^+SO(3)_P\to\Lambda SO(3)_P$ as a
restriction of the
analytic diffeomorphism $(H^-_P)_*\times H_P^+\to (H^-_P)_*\cdot H^+_P$
from Proposition A.1.

Consequently, we have the the induced diffeomorphism
$\tilde\Lambda_*^-SO(3)_P\times\tilde\Lambda^+SO(3)_P\to\tilde\Lambda
SO(3)_P$.
\hfill$\square$

\bibliographystyle{alpha}

Magdalena Toda \\Department of Mathematics and Statistics\\
Texas Tech University \\Lubbock TX 79409
\\e-mail: mtoda@math.ttu.edu

\end{document}